\documentclass[submit,onefignum,onetabnum]{siamonline190516}

\usepackage{lipsum}
\usepackage{amsfonts}
\usepackage{graphicx}
\usepackage{epstopdf}
\usepackage{algorithmic}
\usepackage{harpoon}
\usepackage{caption}
\usepackage{subcaption}
\usepackage{todonotes}
\usepackage{soul} 
\ifpdf
  \DeclareGraphicsExtensions{.eps,.pdf,.png,.jpg}
\else
  \DeclareGraphicsExtensions{.eps}
\fi

\usepackage{enumitem}
\setlist[enumerate]{leftmargin=.5in}
\setlist[itemize]{leftmargin=.5in}

\newsiamremark{remark}{Remark}
\newsiamremark{hypothesis}{Hypothesis}
\crefname{hypothesis}{Hypothesis}{Hypotheses}
\newsiamthm{claim}{Claim}

\headers{Ensemble Inference Methods}{O.R.A. Dunbar, A.B. Duncan, A.M. Stuart, and M.T. Wolfram}

\title{Ensemble Inference Methods for Models With Noisy and Expensive Likelihoods
\thanks{Submitted to the editors DATE.\funding{AD is supported by the UKRI Strategic Priorities Fund under the EPSRC Grant EP/T001569/1,
particularly the ``Digital Twins for Complex Engineering Systems'' theme within that grant, and The Alan Turing Institute.
 MTW is supported by the New Frontier Grant NST-0001 of the Austrian Academy of Sciences.
 AMS is supported by
 NSF (award AGS‐1835860), 
 NSF (award DMS-1818977) and 
 by the Office of Naval Research (award N00014-17-1-2079).
 AMS and MTW are also supported by a
 Royal Society International
 Exchange Grant.}}}

\author{Oliver R. A. Dunbar\thanks{Caltech, 1200 E. California Blvd, Pasadena, CA 91125
  (\email{odunbar@caltech.edu}).}
\and Andrew B. Duncan\thanks{Imperial College London, Huxley Building, South Kensington Campus, SW72AZ London
  (\email{a.duncan@imperial.ac.uk}).}
\and Andrew M. Stuart\thanks{Caltech, 1200 E. California Blvd, Pasadena, CA 91125
  (\email{astuart@caltech.edu}).}
\and Marie-Therese Wolfram \thanks{University of Warwick, Gibbet Hill Road, CV47AL Coventry, UK, 
(\email{m.wolfram@warwick.ac.uk}).}}

\usepackage{amsopn}
\usepackage{color}

\newcommand{\bx}{\overline{\CX}}
\newcommand{\bG}{\overline{\CG}}

\newcommand{\Bz}{\mathcal{B}_0}
\newcommand{\Bo}{\mathcal{B}_1}
\newcommand{\Bt}{\mathcal{B}_2}

\newcommand{\CC}{\mathcal{C}}
\newcommand{\CF}{\mathcal{F}}
\newcommand{\CG}{\mathcal{G}}
\newcommand{\CD}{\mathcal{D}}
\newcommand{\CX}{\mathcal{X}}
\newcommand{\R}{\mathbb{R}}
\newcommand{\Rd}{\mathbb{R}^d}

\newcommand{\Td}{\mathbb{T}^d}
\newcommand{\RNd}{\mathbb{R}^{Nd}}
\newcommand{\TNd}{\mathbb{T}^{Nd}}
\newcommand{\intr}{\int_{\Rd}}
\newcommand{\intt}{\int_{\Td}}
\newcommand{\intrt}{\int_{\Rd \times \Td}}
\newcommand{\eps}{\epsilon}

\newcommand{\intT}{\int_{\mathbb{T}^d}}

\newtheorem{assumption}{Assumption}[section]
\newtheorem{formal}{Formal Perturbation Result}[section]

\definecolor{darkred}{rgb}{.7,0,0}
\definecolor{darkblue}{rgb}{0,0,0.7}
\definecolor{darkgreen}{rgb}{0,0.5,0}
\definecolor{darkpurple}{rgb}{0.5, 0, 0.5}

\newcommand{\ad}[1]{{\color{black}{#1}}}
\newcommand{\mtw}[1]{{\color{black}{#1}}}
\newcommand{\as}[1]{{\color{black}{#1}}}
\newcommand{\od}[1]{{\color{black}{#1}}}
\newcommand{\cG}{\mathcal{G}}
\newcommand{\cN}{\mathcal{N}}

\usepackage{ulem}

\ifpdf
\hypersetup{
  pdftitle={Ensemble Inference Methods for Models With Noisy and Expensive Likelihoods},
  pdfauthor={O. Dunbar, A.B. Duncan, A.M. Stuart, M.T. Wolfram}
}
\fi

\begin{document}

\maketitle

\begin{abstract}

The increasing availability of data presents an opportunity to calibrate unknown parameters which appear in complex models of phenomena in the biomedical, physical and social sciences. However, model complexity often leads to parameter-to-data
maps which are expensive to evaluate and are only available through noisy approximations. This paper is concerned with the use of interacting particle systems for the solution of the resulting inverse problems for parameters. Of particular interest is the case
where the available forward model evaluations are subject to rapid fluctuations, in
parameter space, superimposed on the smoothly varying large scale parametric structure of interest. {A motivating example
from climate science is presented, and ensemble Kalman methods (which \as{do not use the derivative of the parameter-to-data map}) are shown, empirically, to perform well. 
Multiscale analysis is then used to 
analyze the behaviour of interacting particle system algorithms when rapid fluctuations, which we refer to as noise, pollute the large scale parametric dependence of the parameter-to-data map.
Ensemble Kalman methods and Langevin-based methods} (\as{the latter} use the derivative of the parameter-to-data map)
are compared in this light. The ensemble Kalman methods are shown
to behave favourably in the presence of noise in the parameter-to-data map,
whereas Langevin methods are adversely affected. On the other hand,
Langevin methods have the correct equilibrium distribution in the setting of noise-free forward models, 
whilst ensemble Kalman methods only provide an uncontrolled approximation,
except in the linear case. Therefore a new class of algorithms, ensemble
Gaussian process samplers, which combine the benefits of both ensemble Kalman 
and Langevin methods, are introduced and shown to perform favourably.
 
\end{abstract}



\section{Introduction}

The focus of this work is on the solution of inverse problems in the
setting where only noisy {approximations} of the forward problem (the parameter-to-data map)
are available and where the evaluations are expensive.
The methodological approaches we study are all ensemble based.
The take-home message of the paper is that judicious use of 
{ensemble Kalman methodology and generalizations} may be
used to remove the pitfalls associated with gradient based methods
in this setting, but still retain the advantages of gradient
descent; the conclusions apply to both optimization and sampling approaches to inversion.
We provide theoretical and numerical
studies which allow us to differentiate between existing ensemble based approaches, 
and we propose a new ensemble-based method. Subsection \ref{ssec:set} provides the set-up in which we work, Subsection \ref{ssec:lit}
is devoted to a literature review, while Subsection \ref{ssec:ove}
overviews the contributions of the paper and describes
its organization..

\subsection{The Setting}
\label{ssec:set}

The problem we study is this: we seek to infer $x \in \mathbb{R}^d$ given observations $y \in \mathbb{R}^K$ of $G_0(x)$, so that
\begin{align}\label{eq:ob0}
y = G_0(x) + \xi.
\end{align}
The specific instance of the observational noise $\xi$
is not known, but its distribution is; to be concrete we assume that $\xi \sim \mathcal{N}(0, \Gamma)$, with strictly positive-definite covariance $\Gamma \in \mathbb{R}^{K\times K}$. After imposing a prior probability measure
$x \sim \mathcal{N}(m, \Sigma)$, 
application of Bayes rule shows that the resulting
posterior distribution is given by
\footnote{\as{Let $\langle \cdot, \cdot \rangle$, $|\cdot|$ denote Euclidean inner-product and norm. Throughout, for positive-definite symmetric matrix $A$,
we use the notation $\langle \cdot, \cdot \rangle_{A}=
\langle \cdot, A^{-1}\cdot \rangle$, and $|\cdot|_A=|A^{-\frac12}\cdot|.$}}
\begin{align}\label{eq:posteriorb}
    \pi_0(x) &\propto e^{-V_0(x)},\\
    V_0(x) &:=\frac{1}{2}|y - G_0(x)|_{\Gamma}^2+
    \frac{1}{2}|x-m|_{\Sigma}^2.
    \end{align}
This is the standard setting of Bayesian inversion \cite{kaipio2006statistical}. 
The objective is either to generate samples from  target distribution $\pi_0(x)$ (Bayesian approach) or to compute minimizers of $V_0(x)$ (maximum a posteriori estimation -- the optimization approach). The specific focus of this
paper is the setting where $G_0(\cdot)$ is expensive to
evaluate and only a noisy approximation, $G_\eps(\cdot)$,
is available.\footnote{It is important to distinguish between the observational noise $\xi$, appearing in the observations $y$, and
the concept of noisy evaluations of the forward model.} The
parameter $\eps \ll 1$ characterizes the lengthscale, in the space
of the unknown parameter $x$, 
on which the noisy approximation varies.

In order to understand this setting we define 
\begin{align}
\label{eq:forward}
{G}_{\epsilon}(x) = G_0(x) + G_1(x/\epsilon).
\end{align}
Our goal is to solve the inverse problem \eqref{eq:ob0} defined by $G_0$, using only evaluations of $G_\eps$, not of $G_0.$ 
In this context it is also useful to define the multiscale potential
\begin{equation}
\label{eq:potential}
V_\eps(x) :=\frac{1}{2}|y - G_\eps(x)|_{\Gamma}^2+
    \frac{1}{2}|x-m|_{\Sigma}^2,
\end{equation}
and the associated multiscale posterior distribution $\pi_\epsilon \propto \exp(-V_\epsilon)$.
Settings in which $G_1$ is both random and periodic will
be considered. Specifically, we will provide a computational example, arising
in climate modeling, of the setting where $G_1$
\as{represents random fluctuations caused by finite-time
average approximations $G_\eps$ of the desired ergodic averaging
operator $G_0$, demonstrating desirable practical performance of ensemble Kalman methods for both Bayesian and optimization approaches
in this setting. And, in order to provide
deeper theoretical understanding, we will use multiscale
analysis to compare ensemble Kalman algorithms and ensemble 
Langevin algorithms,
for solution of the inverse problem \eqref{eq:ob0} where
$G_1$ is periodic and only evaluation of $G_\eps$ is possible.}

The central message of the paper can now be conveyed by
reference to two different classes of stochastic differential
equations (SDEs), both
defined in terms of $G_\epsilon$, but compared
on the basis of their ability to solve the inverse problem
defined by \eqref{eq:ob0}. The first is the 
ensemble Kalman sampler (EKS) which requires only
evaluations of $G_\epsilon(\cdot)$. The second is the
ensemble Langevin sampler (ELS) which requires evaluations of $V_\epsilon(\cdot)$ and its gradient, and hence requires
evaluations of the action of the gradient of $G_\epsilon(\cdot)$. 
In both cases the ensemble size is $N$. 

The EKS comprises $N$ coupled SDEs in $\R^d$, for $X_t^i$ indexed by
$i=1,\dots,N$, and given by
\begin{equation}
\begin{split}
\label{eq:eksN}
dX_t^i &= -\Bigl(\frac{1}{N}\sum_{n=1}^{N}\langle G_\epsilon(X_t^n) - \overline{G}_{\epsilon,t}, G_\epsilon(X_t^i) - y\rangle_{\Gamma} X_t^n \Bigr)\,dt - C_t\Sigma^{-1}(X_t^i-m) \,dt\\
&\hspace{180pt}\quad\quad\quad+ \frac{d+1}{N}(X_t^i - \overline{X}_t)\, dt 
 + \sqrt{2C_t}\,dW^i_t;
\end{split}
\end{equation}
here the $W^i$ are standard independent Brownian motions in  $\Rd$ and 
\begin{subequations}
\begin{align}
 \overline{X}_t &= \frac{1}{N}\sum_{n=1}^N X_t^n,\qquad
 \overline{G}_{\epsilon,t} = \frac{1}{N}\sum_{n=1}^{N} G_\epsilon(X_t^n),\label{eq:means_intro}\\ C_t &= \frac{1}{N}\sum_{n=1}^{N}\left(X_t^n - \overline{X}_t\right)\otimes \left(X_t^n - \overline{X}_t\right). \label{eq:C_intro}
\end{align} 
\end{subequations}
\as{Thus $\overline{X}_t$ denotes the mean of the ensemble $\{X^i\}_{i=1}^N$, $C_t$ is its empirical covariance} 
and $\overline{G}_{\epsilon,t}$ the mean of the image
of the ensemble under $G_\epsilon.$ 

Using the same notation $X_t^i$ for the ensemble members,
and for
$W_t^i$, independent Brownian motions in $\Rd$, the ELS
may be written as, for $i=1,\dots, N$, 
\begin{equation}
\label{eq:interacting-langevin}
 d{X^i}_t  = -C({X}_t)\nabla {V}_\epsilon({X}^i_t)\,dt + \nabla_{x^i}\cdot C({X}_t)\,dt + \sqrt{2C({X}_t)}dW^i_t.
\end{equation}
Here $C: \R^{Nd} \to \R^{d \times d}$ denotes 
the empirical covariance function of arbitrary
collection of $N$ vectors $\{x^i\}_{i=1}^N$ in $\Rd$
and $X_t=\{X^i_t\}_{i=1}^N$. \footnote{Note that $C(X_t)=C_t$ as defined in $\eqref{eq:C_intro}$; however the
{\em function} $C(\cdot)$ on $\R^{Nd}$ is needed to define
the ELS because of the presence of the divergence contribution
$\nabla_{x^i} C(\cdot)$ in the ELS.}

In the setting where $G_\eps$ is linear the SDEs defining
the EKS and the ELS coincide; \as{they are, however, different from one another in general}.  Ostensibly, both the EKS and ELS as defined above are targeting the distribution $\pi_{\epsilon}$, however they differ drastically in their behaviour as the small length-scale $\epsilon$ goes to zero.  \as{As we shall see, as $\epsilon\rightarrow 0$ the EKS \eqref{eq:eksN} behaves as if $G_\eps$ were
replaced by $G_0$ and hence performs well in recovering
solutions of the inverse problem \eqref{eq:ob0}. In contrast
the ELS \eqref{eq:interacting-langevin} is dominated by the
fluctuations arising from $G_1$ and does not 
perform well in recovering
solutions of the inverse problem \eqref{eq:ob0}.
The EKS effectively denoises $G_\epsilon$ whilst the ELS gets stuck in the noise.} Motivated by the potential success of the EKS for sampling from models with noisy likelihoods, and by the wish to make
controlled approximations of the posterior, \as{we propose here a new class of ensemble method  -- the ensemble Gaussian process sampler (EGPS) -- which can sample effectively from rough distributions without making the ansatz of a Gaussian posterior distribution that is used in the EKS.  The strategy underpinning this method involves evolving an ensemble of particles according to overdamped Langevin dynamics using a surrogate GP emulator to replace the noisy,
and potentially expensive, log-likelihood term.}

\as{The multi-scale analysis and computational experiments that we
present lead to an important dichotomy between different classes of ensemble methods which resonates with the conclusions of \cite{plechavc2019sampling}:  (a) those which calculate the gradient of the log-posterior density for every particle within the ensemble and then aggregate this to update the particle positions; (b) those which evaluate the log-posterior for every particle and then compute a gradient, or approximate gradient.}  
We show that those in class b) are robust to the roughness of the posterior landscape and produce approximations of the posterior \eqref{eq:posteriorb}, using only evaluations
of $G_{\eps}$, but with relaxation times independent of $\eps$;
in contrast the performance of those in class (a) deteriorates as the characteristic length-scale $\eps$ of the roughness converges to zero and do not solve the inverse
problem defined by the smooth component $G_0$, but rather solve a 
different inverse problem exhibiting order one corrections.

\subsection{Literature Review}
\label{ssec:lit}
The focus of this paper is the solution of Bayesian inverse problems, via optimization or probabilistic approaches. Due to the intractability of the posterior distribution associated to a typical Bayesian inversion problem, sampling methods play an important role in exploring the posterior distribution and providing systematic uncertainty quantification.   Due to their wide applicability and practical success, Markov Chain Monte Carlo (MCMC) methods based on Metropolis-Hastings (MH) transition kernels remain the de-facto approach to sampling from high-dimensional and/or complex posterior distributions.  Given sufficient computational effort, an MCMC scheme can return an arbitrarily accurate approximation to an expectation of a quantity of interest, however this often requires large numbers of iterations to provide an accurate characterisation of the posterior distribution \cite{geyer1991markov}.    For Bayesian models with computationally expensive likelihoods, such as those typically arising in  climate modeling \cite{jarvinen2012ensemble}, the geophysical sciences \cite{oliver2008inverse,chandra2019bayeslands} and  agent-based models \cite{gomes2019parameter}, this may render MCMC based methods computationally prohibitive, as they require at least one likelihood evaluation per MCMC step.  

The Ensemble Kalman Sampler \eqref{eq:eksN} is an ensemble-based approach for sampling the posterior distribution associated to a Bayesian inverse problem.  It was introduced in \cite{garbuno2020interacting}, 
without the linear correction term proportional to $d+1$.
The linear correction  was identified
in \cite{nusken2019note,garbuno2020affine} and
ensures that, in the case where the forward map $G$ is linear, the one-particle
marginals of the Gaussian invariant measure deliver the 
solution of the Bayesian linear inverse problem for $G$,
subject to additive noise distributed as $\cN(0,\Gamma)$ and 
subject to prior $\cN(0,\Sigma).$
In the case where $G$ is linear and when initialized with positive-definite initial covariance $C_0$, this
system converges exponentially fast, \as{at problem-independent
rate}, to Gaussian measure
given as solution to the linear inverse problem for $G$
subject to additive noise distributed as $\cN(0,\Gamma)$ and prior $\cN(0,\Sigma)$ 
\cite{garbuno2020interacting,carrillo2019wasserstein}. 
In the nonlinear case, the invariant measure is not known
explicitly; however the output of the finite ensemble SDE may be used as a key component in
other algorithms for solution of the inverse problem \cite{cleary2020calibrate,pavliotis2021derivative}
which come equipped with rigorous bounds for the approximation of the posterior.

When the forward map $G$ is differentiable and its derivative can be computed efficiently, then sampling methods which make use of the gradient of the log-posterior density provide means of exploring the state-space effectively.  For example, one may consider the overdamped Langevin process \cite{pavliotis2015stochastic}, given by the solution of the following SDE:
\begin{equation}
\label{eq:ovlangevin}
    dx_t = -K\nabla V(x_t)\,dt + \sqrt{2K}dW_t.
\end{equation}
Here $K$ is symmetric and positive-definite, but
otherwise is an arbitrary preconditioner.
 Under mild conditions \cite{pavliotis2015stochastic}, the Markov process $(x_t)_{t\geq 0}$, will be ergodic with unique invariant distribution given by $\pi \propto \exp(-V)$, so that $x_t$ will  converge in distribution to $\pi$ as  $t\rightarrow \infty$.   Sampling methods based on discretisations of \eqref{eq:ovlangevin} include the Unadjusted Langevin algorithm (ULA) \cite{roberts1996exponential} \as{as well as its metropolized counterpart, the Metropolis Adjusted Langevin Algorithm (MALA) \cite{bou2010pathwise}, and variants such as the preconditioned Langevin version of the pCN algorithm} \cite{cotter2013} and the Riemmanian Manifold MALA algorithm (RM-MALA)
\cite{girolami2011riemann}. Hybrid (also known as Hamiltonian) Monte Carlo based methods also exploit the gradient of $V$ to explore the state-space \cite{duane1987hybrid}, and have been
generalized to the Riemannian Manifold setting in \cite{girolami2011riemann}. The Ensemble Langevin Sampler \eqref{eq:interacting-langevin}
is defined by allowing
an ensemble of $N$ copies of \eqref{eq:ovlangevin} to interact through a common preconditioner  $C(X_t)$ depending on the solution
of the ensemble of equations.
 Assuming that $C({X}_0)$ is positive definite, then $C({X}_t)$ is positive definite for all $t > 0$ and so ${X}_t$ converges in distribution to $\overline{\pi} = \pi^{\otimes N}$ \cite{garbuno2020affine}. This idea follows from the more general concept of ensemble-based sampling methods which accelerate the Markov chain dynamics by introducing preconditioners computed from ensemble information (e.g. sample covariance) \cite{leimkuhler2018ensemble,garbuno2020affine,carrillo2019wasserstein,duncan2019geometry}.  Since, in the case of a linear forward operator, equation \eqref{eq:eksN} coincides with \eqref{eq:interacting-langevin} \cite{garbuno2020interacting}, this connects the ELS with some
of the previously cited literature on the EKS. Quantitative estimates on the rate of convergence to equilibrium in the setting of preconditioned interacting Langevin equations can be found in \cite{garbuno2020affine,carrillo2019wasserstein}.

A number of other ensemble-based sampling methods have been proposed, building on the Ensemble Kalman Sampler and related work.  In  \cite{pavliotis2021derivative}, the authors propose a multiscale simulation of an interacting particle system, which delivers controllable approximations of the true  posterior;
it is rather slow in its basic form, but can be made more efficient when preconditioned by covariance information derived from the output of the ensemble Kalman 
sampler. Other such ensemble methods include replica exchange  \cite{swendsen1986replica,liu2008monte} as well as  Metropolis-Hastings based approaches \cite{goodman2010ensemble,cappe2004population,jasra2007population,iba2001population}. The recent work \cite{maoutsa2020interacting} also employs an ensemble of particles for evolving density estimation using kernel methods
with the objective of approximating solution of a Fokker-Planck equation.

The presence of multiple modes in the target posterior distribution is a key cause of slow sampler convergence, as any ergodic Markov chain must spend the majority of its time exploring around a single mode, with rare transitions between modes.  Mitigating this issue has  various extensions to standard MH-based MCMC including delayed-rejection methods \cite{green2001delayed,haario2006dram}, adaptive MCMC and methods based on ensembles to promote better state-space exploration, e.g. parallel tempering \cite{liu2008monte} and others.   This issue is further exacerbated for models with posterior distributions exhibiting ``roughness'' characterised by a non-convex, non-smooth posterior with large numbers of local maxima, such as inverse problems 
{arising in climate modeling \cite{dunbar2020calibration}}, Bayesian models in geoscience \cite{chandra2019bayeslands}, frustrated energy landscape models in molecular models of protein structures, glassy models in statistical physics and similar models in the training of neural networks \cite{baity2018comparing}. In the context of Bayesian inverse problems, such pathologies may arise naturally if the forward model exhibits multiscale behaviour \cite{frederick2017numerical}, particularly when only sparse data is available, giving rise to identifiability issues.  Alternatively, this may occur if one only  has access to a noisy estimate of the likelihood, e.g. for some classes of intractable likelihoods such as those arising from stochastic differential equation models with sparse observations.
 Similarly, rough posteriors may also arise if one is fitting a Bayesian inverse problem based on estimators of sufficient statistics of the forward model \cite{morzfeld2018feature}; this setting arises in \as{parameter estimation problems of the type described in \cite{cleary2020calibrate}}, where time-averaged quantities are used for parameter estimation in chaotic dynamical systems. In the special case where one has an unbiased estimator of the likelihood then Pseudo-Marginal MCMC methods \cite{andrieu2009pseudo} provide means of sampling from the exact posterior distribution,  but the performance of these methods degrades very quickly with increasing dimension.  
In the context of uncertainty quantification, Gaussian Processes (GPs) were first used to model ore reserves for mining \cite{krige1951statistical}, leading to the kriging method, which is now well established in the geo-statistics community \cite{stein2012interpolation}. Subsequently, GPs have been very successfully used to provide black-box emulation of computationally expensive codes \cite{sacks1989design}, and in \cite{kennedy2001bayesian} a Bayesian formulation of the underpinning framework is introduced. Emulation methods based on GPs are now widespread, finding applications ranging from computer code calibration \cite{higdon2004combining}, global sensitivity analysis \cite{oakley2004probabilistic}, uncertainty analysis \cite{oakley2002bayesian} and MCMC
 \cite{lan2016emulation}.

Surrogate GP models to accelerate MCMC have been considered before, for example in \cite{lan2016emulation} higher-order derivatives of the log-likelihood required in the calculation of Riemannian Manifold Hamiltonian Monte Carlo were calculated via a GP emulator. Similarly, in \cite{strathmann2015gradient} a nonparametric density estimator based on an infinite dimensional exponential family was used to approximate the log-posterior and then compute the derivatives required for HMC.  Surrogate models to augment existing MCMC methods through a delayed rejection scheme have been considered in \cite{zhang2019accelerating} for GPs and \cite{yan2019adaptive} for neural network surrogates.  In the context of ensemble methods there have been a number of recent works which make use of interpolation in reproducing kernel Hilbert spaces (RKHS) for density estimation and/or gradient estimation which are subsequently used to formulate a ensemble sampling scheme  \cite{reich2019data,maoutsa2020interacting,reich2019fokker,pathiraja2019discrete}.

Our analysis and evaluation of the algorithms
is based on deploying multiscale methodology to determine the
effect of small scale structures on the large scales of
interest; in particular we apply the formal perturbation 
approach to multiscale analysis
which was systematically developed in \cite{bensoussan2011asymptotic}, and
which is presented pedagogically in \cite{pavliotis2008multiscale}.
To simplify the analysis we perform the multiscale
analysis for mean field limit problems, requiring the
study of nonlinear, \as{nonlocal}
Fokker-Planck equations; previous use of
multiscale methods for nonlinear,  \as{nonlocal} Fokker-Planck equations arising from
mean field limits may be found in \cite{gomes2018mean,gomes2019mean}.

\subsection{Our Contributions}
\label{ssec:ove}

In this paper we make the following contributions to the analysis
and development of ensemble based methods for the solution of inverse
problem \eqref{eq:ob0}, based
on forward model $G_0(\cdot)$, given
only access to the noisy approximation $G_{\epsilon}(\cdot)$
in the form \eqref{eq:forward}:

\begin{itemize}

\item we present a parameter estimation problem from climate modeling which naturally leads to the solution of 
an inverse problem of the
form \eqref{eq:posteriorb}, in which only noisy evaluations of
the forward model are available, 
demonstrating favourable behaviour of
ensemble Kalman based methods in a setting where $G_1$ is random;

    \item by means of multiscale analysis in the setting
    where $G_1$ is periodic, we demonstrate that the EKS \eqref{eq:eksN} (which does not use gradients of the forward model) exhibits an averaging property that leads
    to recovery of the SDEs applied with $G_1(\cdot) \equiv 0$;
    
    \item in the same multiscale setting, we demonstrate that the ELS \eqref{eq:interacting-langevin} (which uses gradients of the forward model) exhibits a homogenization property which causes the algorithm to slow down and recover
    an SDE different from the one with $G_1(\cdot) \equiv 0$;
    
    \item we introduce the ensemble GP sampler (EGPS) which combines the benefits of the EKS (averaging out noisy forward model evaluations) with
    the benefits of Langevin based sampling (exact gradients and exact posteriors), overcoming the drawbacks of the two methods
    (uncontrolled approximation of the posterior, slow performance in
    presence of noisy forward model evaluations respectively);
    
    \item we employ numerical experiments to illustrate the averaging and
    homogenization effects of the EKS and Langevin samplers, and to show
    the benefits of the EGPS.

\end{itemize}

The paper is organized as follows.  Section \ref{sec:mot} describes the parameter estimation problem arising in climate modeling that serves to motivate our subsequent
analysis.
In Section \ref{sec:EKS} we define the EKS and study its application to
noisy forward models by means of multiscale averaging. 
In Section \ref{sec:LAN} we define a class of interacting Langevin samplers and study its application to
noisy forward models by means of multiscale homogenization.
Section \ref{sec:EGPS} introduces the new ensemble Gaussian process sampler. Numerical results for all three methods are presented
in Section \ref{sec:NUM} and concluding
remarks are made in Section \ref{sec:CONC}.
\as{The appendices contain details of the climate modeling
example, and of the multiscale analysis.}

\section{Motivating Example}
\label{sec:mot}

\od{We present a specific problem of learning parameters from time-averaged data in an idealized climate model.}
Subsection \ref{ssec:tad} describes the abstract problem
of learning parameters in a dynamical system from time-averaged
data; noisy fluctuations are introduced \od{when} finite
time-averaging \od{is used to approximate a} desired but intractable infinite time-average. 
In Subsection \ref{ssec:gcm} this setting is applied
to the specific motivating example of learning
parameters in a GCM (general circulation model). Subsection \ref{ssec:eka} describes
the application of ensemble Kalman methods to \od{solve} the problem. The example serves
two primary purposes: we provide an explicit instance of a noisy
energy landscape $V_{\eps}$ and we demonstrate the favourable
properties of ensemble Kalman methods when solving optimization
or sampling based inverse problems defined by such a landscape. \od{In particular, the random white noise fluctuations that appear in this example provide a severe test of the ensemble Kalman methodology; in this context the results of this section are very positive regarding the performance of ensemble Kalman methods.}
This motivates the analysis which follows in subsequent sections.

\od{We take great care with our notation in this example. We } label the
unknown $\theta$ and use $\cG_0,\ \cG_1$ and $\cG_\eps$ to denote
the desired forward model, the random fluctuations about it and the
noisy forward model ($\cG_\eps=\cG_0+\cG_1$) to which we have access, respectively. We use \od{the calligraphic} $\cG_\epsilon$, rather than $G_\epsilon$,
to distinguish from the multiscale analysis 
\as{(in subsequent sections)} where the lengthscale $\eps$
is precisely defined through periodic fluctuations in parameter
space . \od{In what follows,} 
$\cG_\eps$ is subject to white noise, and hence
to arbitrarily short lengthscale $\eps$ for the fluctuations;
in practice, however, a value of  $\eps$ satisfying $0<\eps \ll 1$ can still be defined as the minimal separation between the points in $\theta$-space at which $\cG_\eps$ is evaluated. Therefore  $\cG_\eps$ \od{{should be understood with such a choice of $\eps$ }}in the random setting.  

\subsection{Parameter Inference From Time-Averaged Data}
\label{ssec:tad}

Our point-of-departure is the following parameter-dependent dynamical system:
\begin{align}
\label{eq:added}
\frac{du}{ds}&=F(u;{\theta}), \quad u(0)=u_0.
\end{align}
We assume that this dynamical system is ergodic and mixing.
We let $u(t;\theta)$ denote the parameter-dependent solution
of this problem.
Our goal is to learn $\theta$ from data $y$ computed
from finite time-averages of a
\as{function $\varphi(\cdot)$, defined on the state space, over
time-interval of duration $T.$} In detail we have data $y$
given by
\begin{align*}
    y=\cG_{\eps}({\theta})+\xi_{obs},
\end{align*}
where
$$
\cG_{\eps}({\theta})=\frac{1}{T}\int_0^T \varphi\bigl(u(s;\theta)\bigr)ds,
$$
and where $\xi_{obs} \sim \mathcal{N}(0,\Delta_{obs})$ is observational noise.\footnote{It may also be thought of as representing model error, as in \cite{KenOHa01a}.} We note that $\cG_{\eps}$ depends on initial condition $u_0$, which we view as a random variable distributed according to the
invariant measure \as{of \eqref{eq:added}.} For ergodic, mixing dynamical systems a central limit theorem may sometimes be proven to hold \cite{araujo2015rapid}, or empirically observed, for data drawn
at random from the invariant measure. In this setting
\begin{align*}
\cG_{\eps}({\theta})& \approx \cG_0({\theta})+\cG_1({\theta}),\\ \cG_1({\theta}) & \sim \mathcal{N}\bigl(0,T^{-1}\Delta(\theta)\bigr),
\end{align*}
where $\cG_0$ is the infinite time-average which, by ergodicity,
is independent of the initial condition $u_0$. Thus $\cG_\eps(\cdot)$ may be viewed as a noisy perturbation $\cG_1(\cdot)$ (with standard deviation scaling with $T^{-\frac{1}{2}}$) of the function $\cG_0(\cdot)$ \od{where} the noise induced by the unknown initial condition $u_0$ appears only in $\cG_1$ and
not in $\cG_0$. Whenever we evaluate $\cG_\eps$ at different
$\theta$ this noisy evaluation should be thought of as being evaluated
independently with respect to random initial condition $u_0.$ Hence, evaluations of $\cG_\eps$ contain rapid fluctuations that are white in $\theta$-space; as mentioned earlier $\eps$ can be thought of, in
practice, as the minimal separation between $\theta$ values at which time-averages are evaluated.

We approximate $T^{-1}\Delta(\theta)$ by a constant
covariance $\Delta_{model}$ estimated from a single long-run of the (ergodic and mixing) model at a fixed parameter $\theta^\dagger$ and batched into windows of length $T$. If we let $\xi_{model}
\sim  \mathcal{N}(0,\Delta_{model})$, and assume that the
observation error $\xi_{obs}$ is independent of the initial
condition $u_0$, then we have inverse problem
\begin{equation}
\label{eq:gcm_invprob_inf}  
y=\cG_0(\theta)+\xi
\end{equation}
where $\xi=\xi_{obs}+\xi_{model} \sim \mathcal{N}(0,\Delta_{obs}+\Delta_{model}).$
It is then of interest to solve this inverse problem given data $y$
and using algorithms which have access only to (random)
evaluations of $\cG_\eps.$

\subsection{Parameter Estimation in a General Circulation Model}
\label{ssec:gcm}

Climate modeling provides a significant application where the 
set-up of the preceding subsection is relevant. Whilst numerical
weather prediction is entwined with learning the
initial condition $u_0$ of the system, in the setting of climate
modeling $u_0$ is a nuisance parameter of no intrinsic interest.
It is thus natural to calibrate models to time-averaged data,
with goal being the prediction of climate statistics into the
future. In this setting it is natural to solve the inverse
problem relating to infinite time-averages, since the nuisance
parameter $u_0$ disappears from this problem, but to do so given
only the ability to evaluate finite time-averages, because infinite
time-averaging is not feasible in practice.

We consider an idealized moist GCM detailed in \cite{FriHelZur06,OGoSch08}. The GCM comprises three
space-dimensional time-dependent \as{coupled partial
differential equations (PDEs) describing} the large-scale atmospheric motions of an aquaplanet\footnote{A planet covered with a 1m thick slab ocean layer on the inner boundary with no surface topography.}, representing conservation of mass, momentum and
energy. The PDEs are coupled with parameterizations to resolve the subgrid-scale dynamics, notably of moist convection, turbulence and radiation. In the experiments reported, the GCM has a spherical spectral discretization of (32,64,20) discrete latitudes, longitudes, and vertical layers (unevenly spaced in the vertical, with more discrete layers near the planet surface). The time discretization is based on operator-splitting, combining a
leapfrog method (explicit) with a Robert-Asselin time filter  (implicit) -- a standard approach for spectrally discretized atmospheric models \cite{Rob66,Ass72,Wil09} . The simulated aquaplanet is statistically homogeneous in the longitudinal direction, \od{and statistically stationary in time, after an initial burn-in period}. The computational experiments with these discretizations are stable, and take approximately 1-2s per simulated day when distributed over 8 CPU cores.

A Bayesian formulation of the inverse problem of learning two 
subgrid-scale convection parameters, the relative humidity RH and the relaxation time
$\tau$ is presented in \cite{dunbar2020calibration}. The 
quasi-equilibrium moist convection scheme used to describe
subgrid-scale phenomena relaxes temperature and specific humidity towards moist-adiabatic reference profiles with a fixed relative humidity RH.
As data we take three longitudinally averaged (due to symmetry) and time averaged climate statistics: free-tropospheric relative humidity, daily precipitation and the probability of 90th percentile precipitation. These statistics are known to be informative about the unknown parameters \cite{Fri07}. Observing these quantities
latitudinally  \od{and averaging over a window of 30 days,} we obtain data in $\mathbb{R}^{96}$.  The specific formulation of the model, along with details on how the instance of data $y$ and the covariances $\Delta_{obs}$ and $\Delta_{model}$ are constructed are deferred to Appendix \ref{A:details_gcm}.

\subsection{Ensemble Kalman Algorithms For Parameter Inference}
\label{ssec:eka}
We seek to solve the inverse problem \eqref{eq:gcm_invprob_inf} using the ensemble Kalman sampler (EKS) algorithm \cite{garbuno2020interacting}, by evolving a set of $N$ particles $\lbrace \theta^i_{t_n}\rbrace_{i=1}^N$ over a discrete set of 
\as{algorithmic} time steps $0 = t_0 < t_1 < \ldots$\footnote{\as{Not to be confused with physical time $s$ appearing in \eqref{eq:added}.}} such that $\Delta t_n = t_{n+1} - t_n$, to approximate the posterior distribution over $\theta$.   Note that
we distinguish between the EKS, which is an SDE, and the EKS
algorithm, which refers to a discretization of the SDE \eqref{eq:eksN} to
obtain an implementable methodology. 
The mean and covariance under the
prior on the unknown parameter $\theta$ are denoted $(m, \Sigma)$. We define  $M(t_n)=\big(I + \Delta t_n C_{t_n} \Sigma^{-1} \big)$, where $C_{t_n}$ is the empirical covariance of the particles $\lbrace \theta_{t_n}^i\rbrace_{i=1}^N$.  We  consider the approximate posterior sampling approach proposed in \cite{garbuno2020interacting} which has the form,
\begin{subequations}
\label{eq:gcm_eksupdate_dt2}  
\begin{align}\nonumber
  M(t_n)\theta_{t_{n+1}}^{*,i} &= \theta^i_{t_{n}}+ \Delta t_n C_{t_n} \Sigma^{-1} m\\
  &\quad\quad +\Delta t_n\Bigl(\frac{1}{N}\sum_{n=1}^{N}\langle \cG_\epsilon(\theta_{t_n}^n) - \overline{\cG}_{\epsilon,{t_n}}, y - \cG_\epsilon(\theta_{t_n}^i) \rangle_{\Gamma} (\theta_{t_n}^n - \overline{\theta}_{t_n}) \Bigr)  \\
   \theta^i_{t_{n+1}} &=  \theta^{*,i}_{t_{n+1}} + \sqrt{2\Delta t_n C_{t_n}} \xi^i_{n+1},
\end{align}
\end{subequations}
for $i=1,\dots, 20$, where $\xi^i_{n+1} \overset{\text{i.i.d}}{\sim} \mathcal{N}(0,1)$,\footnote{The algorithms we use here do not add additional perturbations of $y$ for each ensemble member, as is often performed \cite{iglesias2013ensemble}.} and  $\Gamma = (\Delta_{model} + \Delta_{obs})I$, and where $\overline{\theta}_{t_n}$ and $\overline{\cG}_{\epsilon,{t_n}}$ are the ensemble averages of $\lbrace \theta_{t_n}\rbrace_{i=1}^N$ and $\lbrace \cG_\epsilon(\theta_{t_n}^i) \rbrace_{i=1}^N$, respectively.  This is a linearly implicit split-step scheme for the SDE \eqref{eq:eks}, but we do not use the finite system size correction,
proportional to \as{$N^{-1}$} and identified in \cite{garbuno2020affine,nusken2019note}, because the effect is small for this example; but it is easily incorporated.   
\as{We also consider the Ensemble Kalman Inversion (EKI) version of this algorithm in
which $M(t_n) \equiv I$ and the white noise contribution is dropped; this corresponds
to time discretization of the ordinary differential equation (ODE)
found by dropping the last three terms on the  righthand side of \eqref{eq:eksN}.} 

\as{We run both the EKS and EKI from algorithmic time $t_0=0$ until (approximately in the EKS case) 
time $5$. For EKI we use $125$ fixed steps, $\Delta t_n=0.04$ for all $n$.} \footnote{\od{In practice, EKI is stable, and produces qualitatively similar solutions over a wide range of choices of $\Delta t_n\leq 1$.}} {For EKS we use an adaptive time-step $\Delta t_n$ that takes values $0.011,0.038,0.11$ for $n=0,1,2$, and $\Delta t_n \approx 0.1$ for $n\geq 3$; we terminate after 49 iterations. For the EKS we use the simple adaptive time-step $\Delta t_n$ proposed  for the Ensemble Kalman Inversion (EKI) algorithm in \cite{KovStu19}, and generalized to the EKS in \cite{garbuno2020interacting}.} We choose $N=20$ and in all our numerical examples we evaluate time-averages with $T=30$ days, the same time-interval used to create the data. An initial ensemble of particles  of size 20 is drawn from the prior.

\begin{figure}[t]
    \centering
    \includegraphics[width=0.45\textwidth]{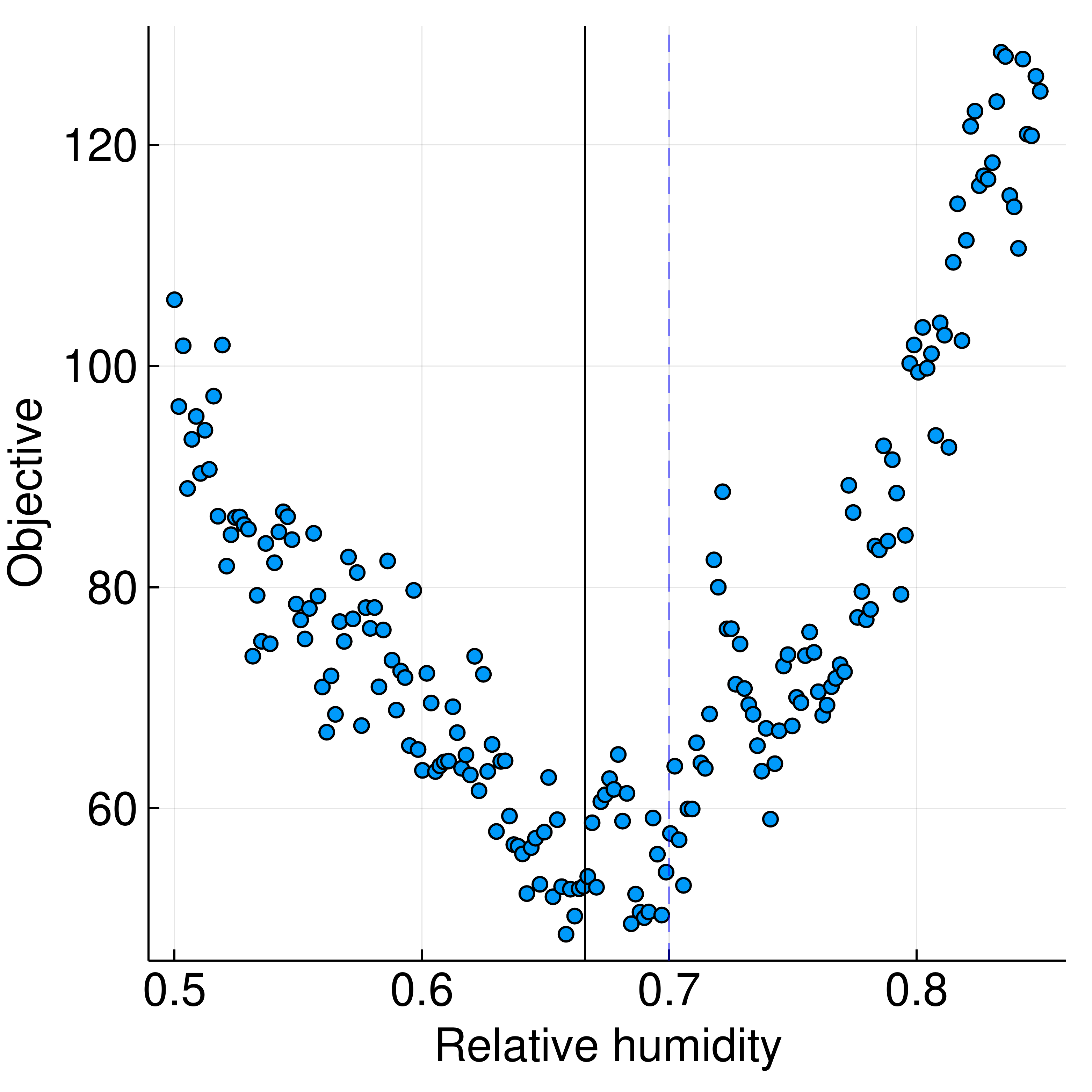}
    \includegraphics[width=0.45\textwidth]{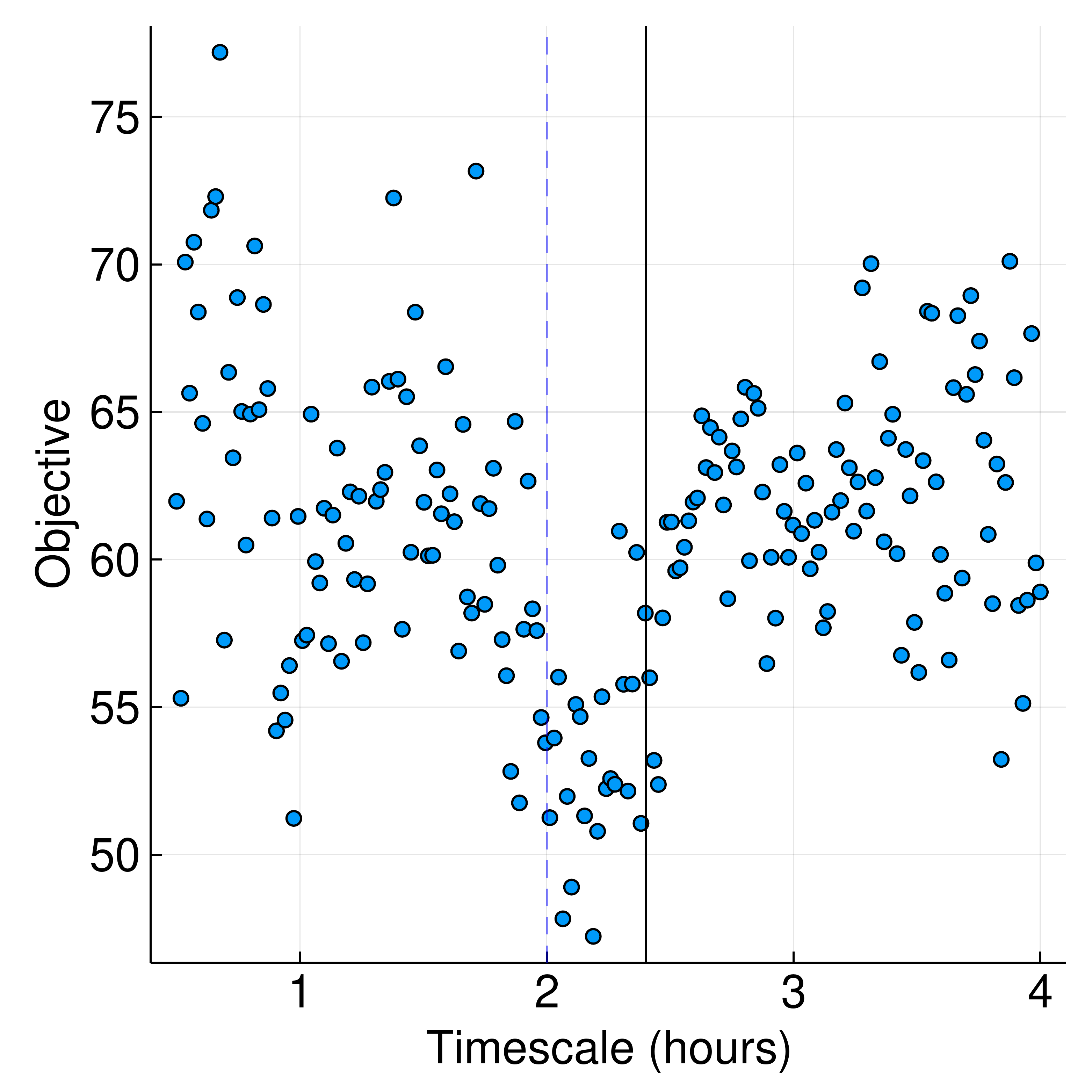}
    \caption{Objective function $V_\eps$ along a line of $200$ parameter values; we vary one parameters and hold the other fixed at the value found from EKI (at time $5$, and for the run $\Delta t=0.04$). The parameter values used to generate the data realization are shown with a blue dashed line; the mean value of the final EKI iteration is a black solid line. \as{The key observation is that EKI produces excellent minimizers despite
    the rough energy landscape.} }
    \label{fig:gcm_fixed_params}
\end{figure}

\begin{figure}
\centering
    \includegraphics[width=0.5\textwidth]{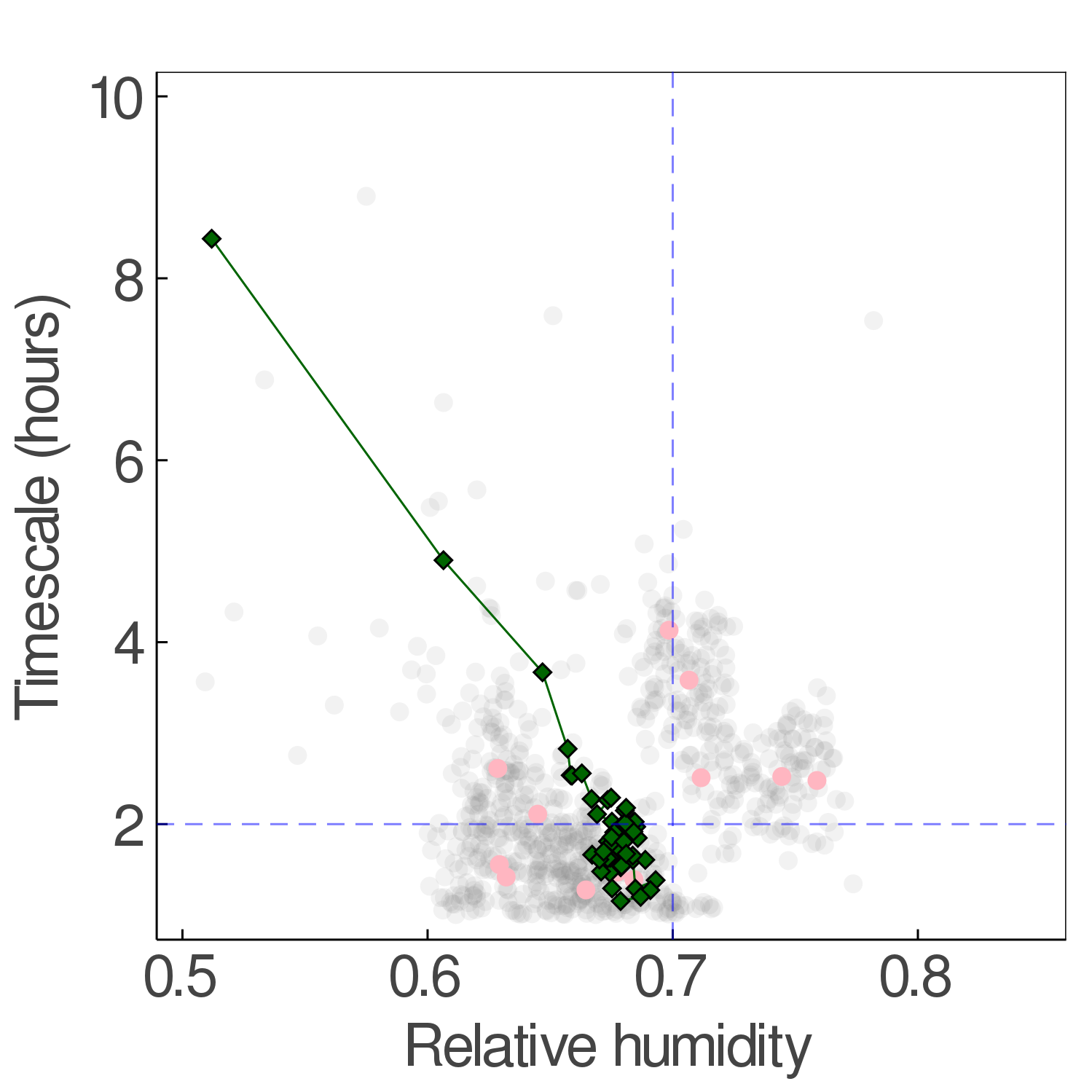}
    \caption{Convergence of the $20$ member EKS ensemble in parameter space over artificial time 0 to 5 with variable time-step \as{as detailed in the main text. The totality of all ensemble members is shown in grey (noting that some points lie outside of the plotting region). The final ensemble (the $49^{th}$)} is given in pink. The parameter set used to generate the data realization is given by the intersection of the blue dashed lines. The green line tracks the ensemble mean over the iterations shown by diamonds. \as{The key observation is that EKS produces excellent samples despite
    the rough energy landscape.}}
    \label{fig:EKS_adaptive}
\end{figure}

We \as{may now} visualize the landscape $V_{\eps}$ 
given by \eqref{eq:potential} with $G_\eps$ replaced by
$G_\eps$. Figure~\ref{fig:gcm_fixed_params}  shows one-dimensional slices through this landscape; in each we hold one parameter at the optimal value of the objective (taken from the EKI
run at time $5$)  while varying the other in uniform increments over $200$ values. The objective evaluations (blue circles) are noisy leading to rapid fluctuations around
a visible convex objective function (defined by the, in practice uncomputable, infinite time-average limit.) Furthermore the optimal parameter set (black vertical line), defined as the mean of the final iteration of the EKI run at time $t =5$, provides a satisfactory 
approximate minimizer of the convex objective function buried underneath the noisy objective function constructed from noisy
forward model evaluations available to the algorithm; and this is
achieved with initialization of the algorithm ensemble from
the prior, which is \od{both broad, and } far from the truth.
 The bias of the objective function with respect to the true parameters (blue dashed vertical line) is to be expected,
 and due to the optimization being performed with a single realization of the noisy data.
 \as{The behaviour of the EKS algorithm is demonstrated in Figure \ref{fig:EKS_adaptive}; it too
 captures the the true parameter very well, despite the noise present in the objective
 function, and also quantifies uncertainty in the estimate, through the spread of the
 pink samples.}

\as{This ability of ensemble Kalman methods to identify approximate minimizers, and generate
approximate posterior samples, in noisy energy landscapes is remarkable.
It leads to the conjecture} that these derivative-free ensemble methods share a common behaviour of ``seeing through the noise'' in model evaluations of ${\cG}_\eps$, enabling
solution of the inverse problem \eqref{eq:gcm_invprob_inf} defined
by ${\cG}_0$. On this basis we will, in the next two sections,
compare derivative-free ensemble methods with gradient-based
interacting particle systems. 

\section{Derivative-Free Sampling}
\label{sec:EKS}

{The previous section demonstrates that algorithms based on
the EKS \eqref{eq:eksN}, and its discretizations  provide remarkable ability to
de-noise rough energy landscapes and identify underlying
smooth landscapes relevant for optimization and sampling, in the
context of a complex problem arising in climate modeling.
In this section we study this problem, returning to the general
set-up of the introduction; we work in continuous time and
with the superimposed rapid fluctuation in the forward model
assumed to be periodic. This simple setting yields
clean theoretical insight,  \od{and experience from the homogenization
and averaging literature \cite{bensoussan2011asymptotic} suggests
that similar results are to be expected for rapid random and periodic
fluctuations, and so can provide explanation of the remarkable behaviour 
observed in the climate modeling example.}
In Subsection \ref{ssec:EKS1} we introduce the mean field limit of \eqref{eq:eksN} which is our starting point; Subsection \ref{ssec:EKS2} is devoted
to analysis of this mean-field SDE, using averaging methods, \as{with detailed
calculations left for an appendix.} The central conclusion from the work in this section is
that the EKS recovers solution to the inverse problem
defined by $\cG_0$ when only evaluation of $\cG_\epsilon$ is available.} 

\subsection{Ensemble Kalman Sampling}
\label{ssec:EKS1}

Here we study the EKS given by \eqref{eq:eksN} where $G_\epsilon$
is defined by \eqref{eq:forward}, and evaluate the relationship
of the SDE to the inverse problem for $G_0$ defined by 
\eqref{eq:ob0}. Our approach is to apply averaging techniques to
the mean field limit of this system.
The mean field limit is given by
\begin{equation}
\label{eq:eks}
dx_t = - \CF(x_t,\rho) \,dt - \CC(\rho)\Sigma^{-1}x_t \,dt  + \sqrt{2\CC(\rho)}\,dW_t,
\end{equation}
where $W$ is a standard Brownian motion \mtw{(independent of the initial condition)} in  $\Rd$ and,
for density $\pi$ on $\Rd$, we define the functions $\overline{\CX},\overline{\CG}$ and
$\CC$ of $\pi$, and $\CF$ of $(\pi,x)$, by
\begin{align}
 \overline{\CX}(\pi) &= \intr X' \pi(X')dX',\quad
 \overline{\CG}(\pi) =  \intr G_{\eps}(X') \pi(X')dX',\\
 \label{eq:cc}\CC(\pi) &= \intr\left(X' - \overline{\CX}(\pi)\right)\otimes \left(X' - \overline{\CX}(\pi)\right) \pi(X')dX',\\
 \CF(x,\pi)&=\Bigl(\intr\langle G_{\eps}(X') - \overline{\CG}(\pi), G_{\eps}(x) - y\rangle_{\Gamma} X'\pi(X')dX'\Bigr).
\end{align} 
Here $\rho$ is the time-dependent density of the process, and self-consistency implies that
it satisfies the nonlinear Fokker-Planck equation
\begin{equation}
    \label{eq:NLFP}
    \partial_t \rho =\nabla_x \cdot \Bigl(\nabla_x \cdot \bigl(\CC(\rho)\rho\bigr) + \CF(x,\rho)\rho\Bigr).
\end{equation}
It is useful to notice that $\CC(\rho)$ depends only on $t$ and not
$x$ and that hence we may also write
\begin{equation}
    \label{eq:NLFPbb}
    \partial_t \rho =\CC(\rho):D^2_x \rho+\nabla_x \cdot \bigl(\CF(x,\rho)\rho\bigr).
\end{equation}
In equation \eqref{eq:NLFP} the outer divergence acts on vectors,
the inner one on matrices; in equation \eqref{eq:NLFPbb} the Frobenius inner-product $:$ between matrices is used.\\
Carrillo and Vaes \cite{carrillo2019wasserstein} established stability estimates in the Wasserstein distance for solutions of \eqref{eq:NLFP} in case of linear $G$, recovering convergence towards equilibrium results by Garbuno et al, see \cite{garbuno2020interacting}.

\subsection{The Small $\epsilon$ Limit}
\label{ssec:EKS2}

In order to understand the performance of the EKS algorithm when rapid fluctuations are present
in the forward model on the EKS algorithm, we proceed to analyze it under the following assumption:

\begin{assumption}
\label{a:1}
The forward model \eqref{eq:forward} satisfies
\begin{align}\label{eq:aforward}
    {G}_{\epsilon}(x) = G_0(x) + G_1(x/\epsilon),
\end{align}
$G_0 \in C^{1}(\mathbb{R}^d,\mathbb{R}^K)$,  $G_1
\in C^{1}(\mathbb{T}^d,\mathbb{R}^K)$ and $\int_{\mathbb{T}^d} G_1(y)\,dy = 0$. $\quad\qed$ 
\end{assumption}

Here $\Td$ denotes the $d$ dimensional unit torus:  $G_1$ is a $1-$periodic function in every direction. Although the periodic
perturbation is a simplification of the typical noisy models
encountered in practice, \as{such as the class presented in Section \ref{sec:mot}}, it provides a convenient form for analysis which is enlightening about the behaviour of algorithms more generally; furthermore the multiscale ideas we use may be generalized to stationary
random perturbations and similar conclusions are to be expected
\cite{bensoussan2011asymptotic}.

We use formal multiscale perturbation expansions to understand the
effect of the rapidly varying perturbation $G_1(\cdot)$ on the
smooth envelope of the forward model, $G_0(\cdot)$, in the context of the
EKS, using the mean field limit.
To describe the result of this multiscale analysis we define the
averaged mean field limit equations, found from \eqref{eq:eks}
and \eqref{eq:NLFP} with $G_1(\cdot) \equiv 0$ so that
$G_\eps(\cdot)$ may be replaced with $G_0(\cdot)$:
\begin{equation}
\label{eq:eks0}
dx_t = - \CF_0(x_t,\rho_0) \,dt - \CC(\rho_0)\Sigma^{-1}x_t \,dt  + \sqrt{2\CC(\rho_0)}\,dW_t,
\end{equation}
with
\begin{align*}
 \overline{\CG}_0(\pi) &=  \intr G_0(X') \pi(X')dX',\\
 \CF_0(x,\pi)&=\intr\langle G_0(X') - \overline{\CG_0}(\pi), G_0(x) - y\rangle_{\Gamma} X'\pi(X')dX'.
\end{align*} 
To be self-consistent the density $\rho_0(x,t) \in C\bigl((0,\infty);L^1(\Rd;\R^+)\bigr)$ must satisfy 
the nonlinear Fokker-Planck equation
\begin{equation}
    \label{eq:NLFP0}
    \partial_t \rho_0 =\nabla_x \cdot \Bigl(\nabla_x \cdot \bigl(\CC(\rho_0)\rho_0\bigr)+\CF_0(x,\rho_0)\rho_0\Bigr).
\end{equation}
The following result is derived in Appendix \ref{A:1} and it
shows that, as $\epsilon \to 0$, equation \eqref{eq:NLFP} is approximated
by \eqref{eq:NLFP0}.

\begin{formal} \label{fpr:1}
Let Assumption \ref{a:1} hold with $0< \eps \ll 1$. If the solution of \eqref{eq:NLFP} is expanded in the form $\rho=\rho_0+\eps \rho_1+\eps^2 \rho_2+\cdots$, then formal multiscale analysis demonstrates
that $\rho_0$ satisfies \eqref{eq:NLFP0}.
\end{formal}

\begin{remark}
\label{rem:1}

\begin{itemize}
\item The result shows that, as $\epsilon \to 0$, the mean field SDE \eqref{eq:eks}, and the nonlinear Fokker-Planck equation \eqref{eq:NLFP} for its density, are approximated by the SDE \eqref{eq:eks0}, and the nonlinear Fokker-Planck equation \eqref{eq:NLFP0} for its density.
This means that the EKS algorithm simply behaves as if $G_1 \equiv 0$, and ignores the rapid ${\mathcal O}(1)$ fluctuations on
top of $G_0$; this is a very desirable feature for computations whose goal is to solve the inverse
problem \eqref{eq:ob0} defined by $G_0$ but where only black box evaluations of $G_{\epsilon}$ given by
\eqref{eq:forward} are available.

\item \as{This result is consistent with what we observed 
empirically in the behaviour of ensemble Kalman based algorithms used to learn parameters in a GCM.}

\item We choose to formulate this result in terms of the mean field limit because this leads to a
transparent derivation of the relevant result. The analysis is cleaner in this limit
as it concerns a nonlinear Fokker-Planck equation with spatial
domain $\Rd \times \Td$; similar results may also be obtained
for the finite particle system by considering a linear 
Fokker-Planck equation with spatial domain $\RNd \times \TNd$.

\item Rigorous justification of the formal expansion could be approached
by \mtw{using the} It\^o formula (see Chapters 17 and 18 in \cite{pavliotis2008multiscale} for example); the main technical
difficulty in this setting is the need to derive bounds
from below on the covariance operator, something which is
considered in \cite{garbuno2020affine} where the finite particle
system is proved to be ergodic.

\end{itemize}
\end{remark}

\section{Derivative-Based Sampling}
\label{sec:LAN}

\as{We now study the ELS \eqref{eq:interacting-langevin}
and study its relation to solution of the Bayesian inverse problem
defined by  \eqref{eq:ob0}.
In Subsection \ref{ssec:ELS} we
introduce the mean-field limit of the ELS, which is our starting point; Subsection \ref{ssec:ELS2} is devoted
to analysis of this mean-field SDE, using homogenization methods
with detailed
calculations left for an appendix.
The central conclusion from the work in this section is that, in contrast to the EKS studied in the
previous section, the ELS performs poorly at recovering
solution to the inverse problem defined by $G_0$ when only
evaluation of $G_\eps$ is available.}

\subsection{Ensemble Langevin Sampling}
\label{ssec:ELS}

The mean field limit of the ELS \eqref{eq:interacting-langevin}
takes the form
$$
	dx_t = -\mathcal{C}(\rho_t)\nabla V_\epsilon(x_t) + \sqrt{2\mathcal{C}(\rho)}dW_t,
$$
where  function $\mathcal{C}(\cdot)$ on densities is as in \eqref{eq:cc} and $V_\epsilon$ is given in \eqref{eq:potential}.
By self-consistency, the associated non-linear Fokker-Planck equation for the time-dependent density of the process $\rho \in C((0, \infty); L^1(\mathbb{R}^d ; \mathbb{R}^+))$ is given by
\begin{equation}
\label{eq:ensemble_langevin}
\partial_t \rho = \nabla_x \cdot\bigl(\mathcal{C}(\rho)\left(\nabla_x\rho+\nabla_x V_\epsilon\rho\right)\big).
\end{equation}
Similarly to the previous section, this may be rewritten as
\begin{equation}
    \label{eq:NLFPbbb}
    \partial_t \rho =\CC(\rho):\Bigl(D^2_x \rho+\nabla_x \bigl(\nabla_x V_\epsilon\rho\bigr)\Bigr).
\end{equation}

\subsection{The Small $\epsilon$ Limit}
\label{ssec:ELS2}

As an ensemble scheme, the system described by \eqref{eq:interacting-langevin} aggregates information from individual particles to obtain a better informed direction in which to explore the posterior distribution.   Unlike the EKS, these approaches compute the gradient before aggregating across particles.  We show that this causes the resulting sampler to be poorly performed with respect to the
presence of rapid fluctuations in the evaluation of the likelihood. 
The following result is derived in Appendix \ref{A:2}. It characterises the evolution of the $O(1)$ leading order term  of $\rho$ solving 
\eqref{eq:ensemble_langevin}. Unlike the
setting in the previous section for
the EKS, the limit is not the same as the Fokker-Planck equation obtained 
from applying the ELS methodology to the inverse
problem defined by forward model
$G_0$ with posterior given by
\eqref{eq:posteriorb}.

\begin{formal} \label{fpr:2}
Let Assumption \ref{a:1} hold with $0< \eps \ll 1$. If the solution of \eqref{eq:ensemble_langevin} 
is expanded in the form $\rho=\rho_0+\eps \rho_1+\eps^2 \rho_2+\cdots$, then formal multiscale analysis demonstrates
that $\rho_0$ satisfies the following mean field PDE:
\begin{equation}
\label{eq:effective_mfe}
\partial_t \rho_0 = \nabla_x\cdot\left(\mathcal{D}(\rho_0)\left(\nabla_x
\rho_0 + \nabla_x \overline{V}\rho_0\right)\right),
\end{equation}
where $\overline{V} =  V_0 - \log Z(x) $,  $Z(x) = \int_{\mathbb{T}^d} e^{-V_1(x,z)}\,dz$ 
and 
$$
\mathcal{D}(\rho_0) =\frac{1}{Z(x)}\int_{\mathbb{T}^d}\bigl(I+\nabla_z\chi(x,z)\bigr)^{\top} \mathcal{C}(\rho_0)\bigl(I+\nabla_z\chi(x,z)\bigr)e^{-V}\,dz.
$$
Here $\chi:\mathbb{R}^d\times \mathbb{T}^d \rightarrow \mathbb{R}^d$ is a solution to the following second order PDE in $z$, parameterized by $x$:
$$
\nabla_z\cdot\left(\mathcal{C}(\rho_0)e^{-V(x)}(\nabla_z \chi(x,z) + I) \right) = 0, \quad (x, z) \in \mathbb{R}^d\times\mathbb{T}^d.
$$
Furthermore,
for arbitrary $\zeta \in \mathbb{R}^d$,
\begin{equation}\label{eq:oqf}
 \zeta^\top \mathcal{D}(\rho_0)\zeta \leq  \zeta^\top  \mathcal{M}(\rho_0)\zeta.
\end{equation}
\end{formal}
\begin{remark}
\begin{itemize}
    \item The homogenized mean field equations in the $\epsilon\rightarrow 0$ limit describe the evolution of a density $\rho_0$ with unique invariant distribution given by $\overline{\pi}(x)\propto \pi_0(x) Z(x)$.  This invariant distribution will generally not be equal to the invariant distribution $\pi_0$, associated with the smoothed inverse problem \eqref{eq:ob0}, defined in \eqref{eq:posteriorb}
    because of the presence of $Z(x)$. This indicates that using an ensemble of coupled Langevin particles applied with potential $V_\epsilon$ derived from the noisy forward problem $G_\epsilon$ will not result in an `averaging out' to obtain samples from the posterior of the smoothed inverse problem with potential $V_0$ derived from the smooth forward problem $G_0$; indeed there will in general be  an $O(1)$ deviation from the target invariant distribution.  
    
    \item A second effect that is caused by the fast-scale perturbation is a slow-down of convergence to equilibrium, specifically  \eqref{eq:oqf} implies that the spectral gap associated with the mean field equation \eqref{eq:effective_mfe} will be generally smaller than that associated with the slowly-varying forward operator $G_0$.
    
    \item The same considerations described in the \as{third and fourth bullets} of Remark \ref{rem:1} also apply here.
    
    \end{itemize}
\end{remark}

\section{The Best Of Both Worlds}
\label{sec:EGPS}
\as{In this section we detail a  gradient-free ensemble method which makes use of smooth estimates of the log-likelihood over the ensemble of particles to estimate the gradient from the available
noisy log-likelihood evaluations.  This approximation is then used to evolve each particle forward according to overdamped Langevin dynamics in the implied smoothed potential. The proposed method has the advantage of
the EKS (robust to noisy energy landscapes) and of the ELS
(works with gradients and provides controllable approximation
of the invariant measure). } \ad{In particular, we expect the convergence to the invariant distribution to be faster compared to EKS to due the exploitation of the approximate gradient, which is insensitive to the local noise-induced fluctuations.}

To this end, we model the partially observed potential
$$V_L(x) =\frac{1}{2}\langle y - G(x), \Gamma^{-1}(y -G(x))\rangle$$ 
as a Gaussian process $f \sim GP(0,k)$,
where $k$ is an appropriately chosen positive definite kernel on $\mathbb{R}^d$. 
This idea is inspired by the paper \cite{maoutsa2020interacting} which uses a closely related approach, with the
goal of approximating solutions to a Fokker-Planck equation.
In this work, we choose $k$ to be a Gaussian radial basis function kernel of the form $k(x,y; \lambda, l)= \lambda \exp(-\lVert x - y\rVert^2/2l^2)$, where $\lambda > 0$ is the kernel amplitude and $l > 0$  is the kernel bandwidth. 
\as{Given (noisy) evaluations of the potential at the ensemble of points} ${X}_t = (X_t^1, \ldots, X_t^N) \in \mathbb{R}^{N\times d}$ we seek a function $f$ such that, \as{for
some $\sigma>0$,}
$$
V_L(X^i_t) = f(X^i_t) + \sigma\xi^i, \quad \xi=(\xi^1,\cdots,\xi^N) \sim \mathcal{N}(0, I).
$$
The corresponding Gaussian process posterior  for $f$ has mean function
$$\widehat{V_L}\textbf{}(x;  \sigma, \lambda, l) = \sum_{i,j=1}^N k(x, X^i_t; \lambda,  l)K({X}; \sigma,\lambda, l)^{-1}_{ij}V_L(X^j_t), \quad x\in\mathbb{R}^d$$ 
and covariance function
$$
\gamma(x,y;  \sigma, \lambda, l) = K(x,y;  \sigma, \lambda, l)-\sum_{i,j=1}^N k(x, X^i_t; \lambda, l)K({X};   \sigma, \lambda, l)^{-1}_{ij} k(X^j_t, y;  \lambda, l).
$$
Here $K({X})_{i,j} = \sigma^2 \delta_{i,j} + k(X^i_t, X^j_t)$.   The gradient of the posterior mean is well-defined and given by
$$
   \nabla\widehat{V_L}(x;   \sigma, \lambda, l) =  \sum_{i,j=1}^N\nabla_x k(x, X^i_t; \lambda, l)K({X};  \sigma, \lambda, l)^{-1}_{ij}V_L(X^j_t).
$$
The particles in the ensemble are then evolved forward according to overdamped Langevin dynamics, i.e.
\begin{align}\label{eq:gpf}
dX_t^i =  -\nabla\widehat{ V_L}\textbf{}(X_t^i;   \sigma, \lambda, l)\,dt - \Sigma^{-1}X_t^i\,dt + \sqrt{2}\,dW_t.
\end{align}

\as{In simply situations the learned energy $\widehat{ V_L}$ is updated every time-step. The three hyperparameters $(\sigma, \lambda, l)$ are chosen to reflect the spread and local variation in the data and hence,
as the conditioning points are updated, these parameters are also adjusted accordingly.} We impose log-normal priors on the amplitude $\lambda$ and observation noise standard deviation $\sigma$ and a Gamma prior on the lengthscale $l$. These prior modeling choices on the hyper-parameters ensure that the posterior mean does not introduce any short-scale variations below the levels of the available data  \cite{cleary2020calibrate,gelman2013bayesian,gelman2017prior}. As is standard in the training of GPs we center
and rescale the training data to have mean zero and variance one. To select 
the hyperparameters we employ an empirical Bayesian approach: we compute the maximum a-posteriori values of the
hyperparameters after marginalising out $f$. This entails selecting $(\sigma, \lambda, l)$ which maximise the log marginal posterior,
\begin{align*}
MLP(\sigma, \lambda, l; {X}) &\propto   \frac{1}{2}\log \sum_{i,j=1}^N {\widehat{V_L}}(X^i_t;  \sigma, \lambda, l) K({X};  \sigma, \lambda, l)^{-1}{\widehat{V_L}}(X^j_t;  \sigma, \lambda, l) \\
&\quad - \frac{1}{2}\log\det K({X};  \sigma, \lambda, l) + \log p_0(\sigma, \lambda, l),
\end{align*}
where $p_0$ denotes the prior density over the hyperparameters.
\\\\
In simulations we employ an Euler-Maruyama discretisation of \eqref{eq:gpf}, coupled with a gradient descent scheme for adaptively selecting the hyperparameters. Let $X_n = (X_n^1, \ldots, X_n^N) \in \mathbb{R}^{N\times d}$ denote the particle ensemble at time-step $n$.  The algorithm for evolving the particles forward to time-step $n+1$ is summarised as follows:
\begin{itemize}
    \item For $i=1,\ldots, N$:
        \begin{itemize}
            \item Set $X_{n+1}^i =X_n^i -  \Delta t\nabla \widehat{V_L}\textbf{}(X_n^i;  \sigma_n, \lambda_n, l_n)\, - \Delta t \Sigma^{-1}X_n^i\,  + \sqrt{2\Delta t}\,\xi_{n},$ where $\xi \sim \mathcal{N}(0,1)$ iid.
            \end{itemize}

            \item Update $(\sigma_{n+1}, \lambda_{n+1}, l_{n+1}) = (\sigma_{n}, \lambda_{n}, l_{n}) + \delta t \nabla_{(\sigma,\lambda,l)}MLP(\sigma_n, \lambda_n, l_n; X_{n+1})$.
\end{itemize}
In the above $\Delta t$ and $\delta t$ are step-sizes for the Langevin updates and the hyperparameter gradient descent, respectively.  \ad{The choice of $\Delta t$ and $\delta t$ is problem-dependent. While it is possible it is possible that  integration with EGP (and EKS) is stiff during the initial transient phase, this can be remedied by using simple adaptive time-stepping. Moreover, for the EGP we see that one of the effects of the smoothing effect of the kernel
is to reduce scale separation within the posterior, thus resulting in less stiff dynamics.
}  

If we are in a situation where evaluating the likelihood is computationally intensive, then we may consider a straightforward modification of these dynamics where time is split into a fixed set of epochs where we keep the same conditioning points within the same epoch, performing several steps of Langevin updates and hyperparameter tuning based on the same conditioning points.   This permits effective exploration of the posterior distribution but with a fixed number of log-likelihood evaluations. \ad{Note that, while we have based the proposed scheme on an underlying Gaussian process model, the use of other non-parametric regression models would be possible, provided that gradients can be readily computed.}

\section{Numerical Results}
\label{sec:NUM}
\as{In this section
we illustrate the performance of the different methods analyzed
or introduced in the preceding three sections, comparing
their performance on three different numerical examples. In particular we compare the EKS defined in equation \eqref{eq:eksN},  the ELS as defined in equation \eqref{eq:interacting-langevin} and the  EGPS as introduced in Section \ref{sec:EGPS}. Our results show the desirable behaviour of the EKS with respect to its ability to avoid the rapid fluctuations imposed
on the smooth parametric structure of interest and rapidly converge to the desired smooth posterior; they show the undesirable slowing down of the ELS, but do not illustrate the
modified limiting posterior as their slow performance means that this equilibrium distribution is not reached in a reasonable number of iterations. They also demonstrate that the EGPS has the same quality of performance as the EKS, with  further improved rate of convergence in continuous time; however, in the
units of evaluations of the (assumed expensive) forward model the EKS
may still remain competitive. The three examples considered are a
perturbed linear model, in Subsection \ref{ssec:plm},
the Lorenz '63 system with parameter estimation
through time-averaged data (as for the GCM) in
Subsection \ref{ssec:lorenz} and a multi-modal
example, in Subsection \ref{ssec:multimodal}.}

\subsection{A Linear Model}
\label{ssec:plm}
As a first pedagogical example we consider solving the inverse problem of the form \eqref{eq:ob0} for a forward map $G_{\epsilon}$ of the form,
for $x=(x_1,x_2)$,
\begin{align}
{G}_{\epsilon}(x) & = {G}_{0}(x)+{G}_{1}(x/\epsilon),\nonumber\\
    {G}_{0}(x) &= A x, \quad
    {G}_{1}(x) =  \left[\sin \left(2 \pi x_1\right), \sin \left(2 \pi x_2 \right)\right]^\top,\label{eq:fwd_linear}\\
     A &= \begin{pmatrix} -1 & 0 \\ 0 & 2 \end{pmatrix}.\nonumber
\end{align}
The objective is to recover the posterior distribution associated with  ``slowly-varying'' component of the forward model $G_0(x) = Ax$,  based on evaluations of the multiscale forward map $G_{\epsilon}$.  To this end, we generate observed data $y \in \mathbb{R}^2$ for a true value of $x$ given by $x^{\dagger}=(-1,1)$ and observational covariance $\Gamma = \gamma^2 I$ with   $\gamma^2 = 0.05$. We impose a Gaussian prior on the unknown parameter $x$ with zero mean and covariance $\Sigma=\sigma^2 I$, with $\sigma^2 = 0.05$. \as{In the absence of the multiscale perturbation $G_1$ in the forward model, the resulting posterior is Gaussian with mean and covariance  given by}
$$
    m_{post} = \left(\begin{matrix}\frac{-1}{1+\sigma^2} & 0 \\ 0 & \frac{2}{4+\sigma^2}\end{matrix}\right)y,\quad \mbox{ and }\quad C_{post} = \gamma^2\left(\begin{matrix}\frac{1}{1+\sigma^2} & 0 \\ 0 & \frac{1}{4+\sigma^2}\end{matrix}\right),
$$
respectively.  Setting $\epsilon = 0.1$, each of the the \as{three} methods is used to evolve an ensemble of $1000$ particles from a $U[0,1]^2$ initial distribution, over a total of $10$ time-units.  The step-size employed for each method is selected differently to ensure stability. 

\begin{figure}
    \centering
    \begin{subfigure}[b]{0.45\textwidth}
    \includegraphics[width=\textwidth]{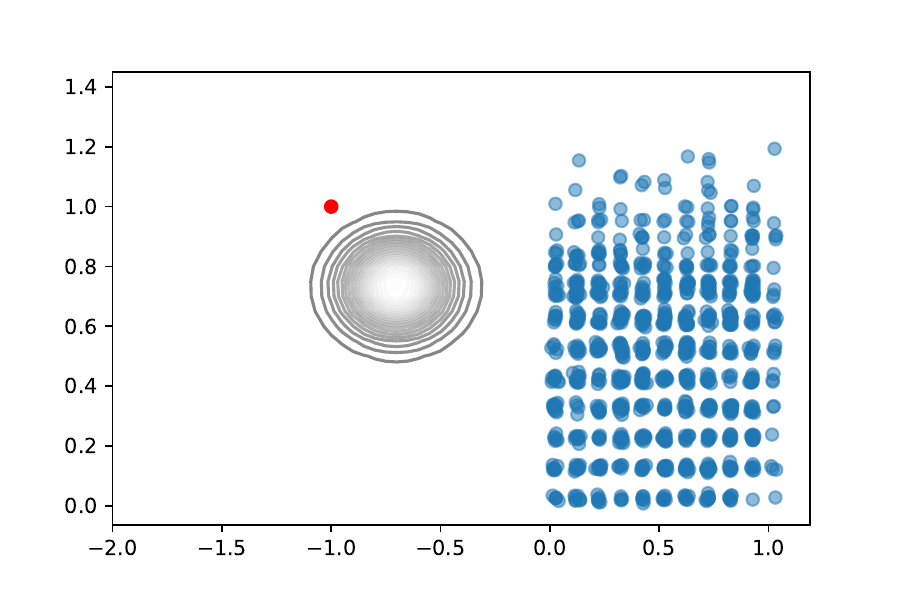}
    \caption{Ensemble Langevin Sampler}\label{fig:lin_ELS}
    \end{subfigure}
    \hfill
    \begin{subfigure}[b]{0.45\textwidth}
    \includegraphics[width=\textwidth]{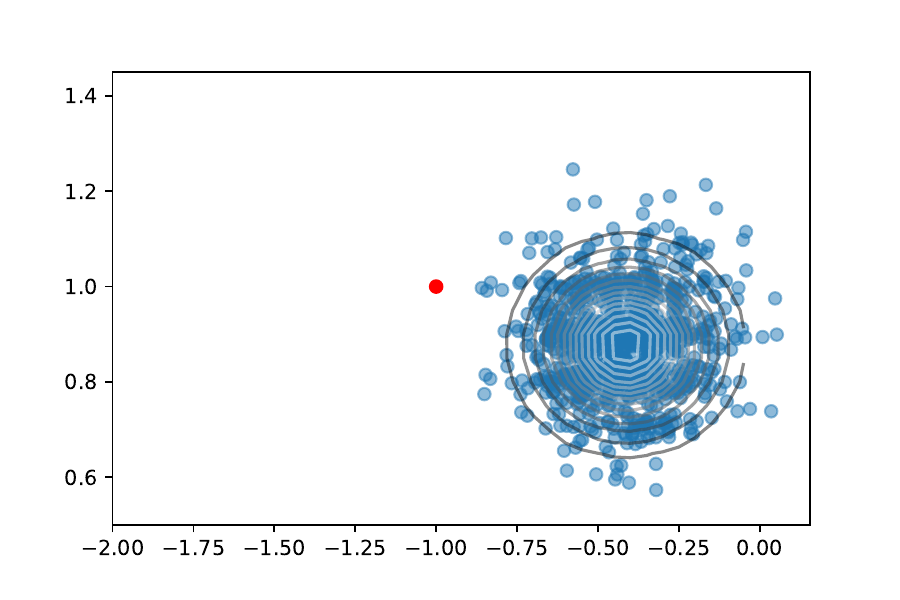}
    \caption{Ensemble Kalman Sampler}\label{fig:lin_EKS}
    \end{subfigure}
    
    \begin{subfigure}[b]{0.45\textwidth}
    \includegraphics[width=\textwidth]{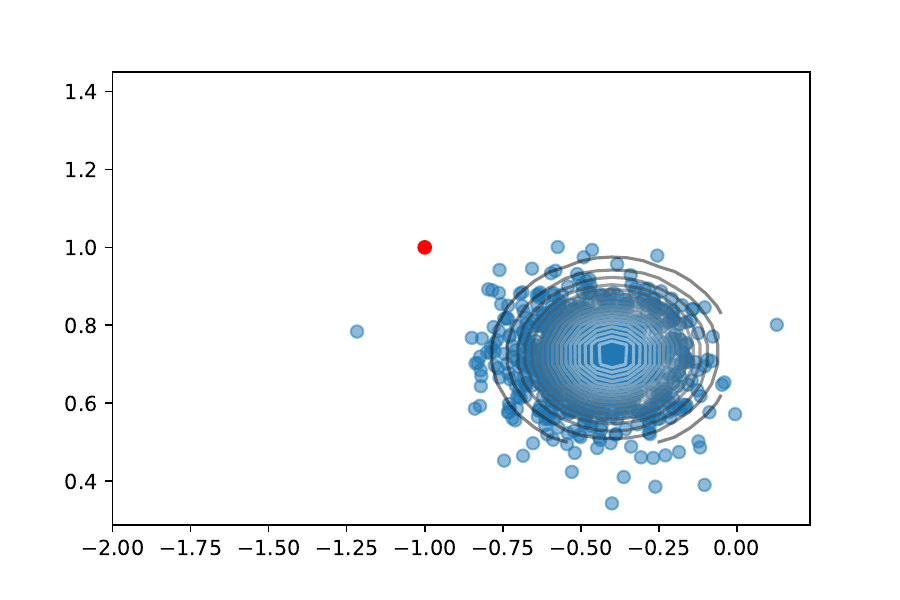}
    \caption{Ensemble GP Sampler} \label{fig:lin_GP}
    \end{subfigure}
    \hfill
    \begin{subfigure}[b]{0.45\textwidth}
    \includegraphics[width=\textwidth]{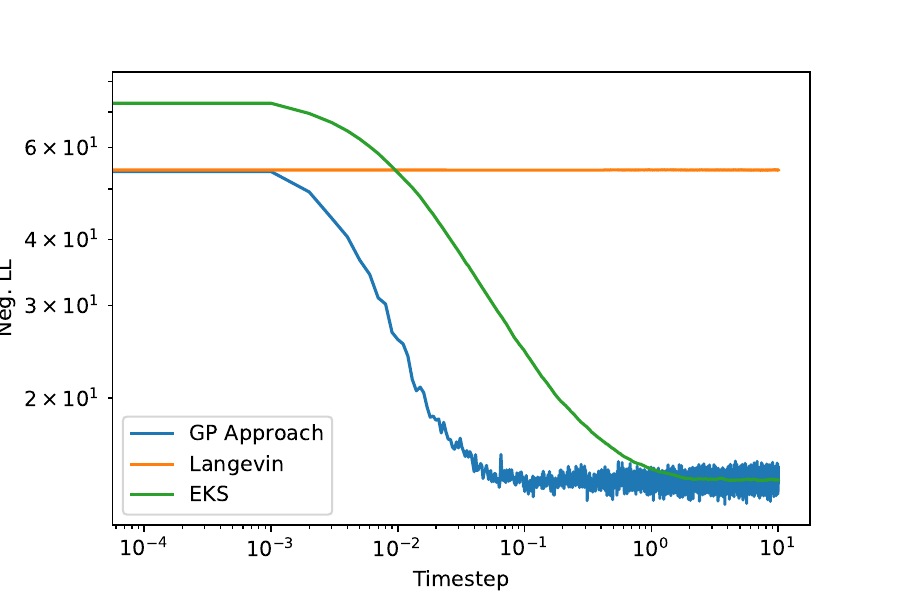}
    \caption{Negative Log-Likelihood} \label{fig:lin_NLL}
    \end{subfigure}
    
    \caption{Plot of the particle ensemble for each of the three processes after simulating $10$ time units.  The contour plot indicates the posterior distribution for the ``slowly-varying'' forward model.  The red dot denotes the truth. Note the different scaling of the axes in plots (a-c). The bottom left plot shows the evolution of negative log likelihood as a function of time-step.}
    \label{fig:linear_hist}
\end{figure}

Figure \ref{fig:linear_hist} shows the particles ensemble at the final time (blue), the true solution (red dot) as well as the posterior $\pi_0$ \eqref{eq:posteriorb} associated with the slowly varying part of the forward model $G_0$ (black contour lines). We observe that particles get stuck in the many local minima \as{for the ELS, shown in Figure \ref{fig:lin_ELS}. This is consistent with the Formal
Perturbation Result \ref{fpr:2}, which indicates that the multiscale perturbation will slow down the dynamics and will result in a significant deviation of the invariant measure of the SDE from the case $G_1=0$.} This is not the case for both the EKS and EGPS, which are able to recover the slowly varying target distribution correctly, see Figure \ref{fig:lin_EKS} and Figure \ref{fig:lin_GP} respectively. Figure \ref{fig:lin_NLL} shows the negative log-likelihood for $\mathcal{N}(m_{post}, C_{post})$ averaged across all the particles in the ensemble,
\as{as the algorithm progresses. Both the EKS and EGPS rapidly move towards the mode of the slowly-varying target distribution, with the EGPS 
converging faster due to the use of the approximate gradient.  On the other hand the Langevin sampler is strongly influenced by the multiscale perturbations, and after $10$ time units remains distant from the desired Gaussian posterior distribution, stuck in
local minima caused by fluctuations $G_1.$ While
the EGPS converges the fastest, it is sensitive to the initial selection of hyper-parameter values and step-sizes; these were initially set through a preliminary tuning phase for the displayed results.  On the other hand, the EKS is remarkably robust to the choice of step-size.} 

\subsection{The Lorenz-63' model}
\label{ssec:lorenz}

Here we work in the setting of parameter inference from time-averaged data, introduced in Subsection \ref{ssec:tad}. Rather than work with the complex
GCM studied in Section \ref{sec:mot}, we work
with the Lorenz '63 model in order to present
controlled experiments at cheaper cost.

The 3-dimensional Lorenz equations \cite{lorenz1963deterministic} are given in the form
\begin{subequations}\label{eq:lorenz63}
\begin{align}
    \dot{x}_1 &= \sigma (x_2 - x_1)\\
    \dot{x}_2 &= r x_1 - x_2 - x_1 x_3\\
    \dot{x}_3 &=x_1 x_2 - b x_3,
\end{align}
\end{subequations}
with parameters $\sigma,\, r, \, b \in \mathbb{R}_+$. In the following we fix the parameter $\sigma=10$ and focus on the inverse problem of identifying $r$ and $b$ from time-averaged data. To this end, we impose a multivariate log-normal prior on $\theta=(r,b)$ with mean $m = (3.3, 1.2)$ \as{and covariance $\Sigma = \text{diag}(0.15^2,0.5^2)$; to be concrete this defines the prior distribution satisfied by $\log(\theta)$.}

  In the notation of Subsection \ref{ssec:tad} we take $T=10$ and define $\varphi \colon \mathbb{R}^3 \rightarrow \mathbb{R}^9$ given by
\begin{align*}
    \varphi(x) = (x_1, x_2, x_3, x_1^2, x_2^2, x_3^2, x_1 x_2, x_2 x_3, x_1 x_3);
\end{align*}
this defines forward model $\cG_\eps$ as a time-average of first and
second moments of the solution over $10$ time units. 
In the experiments that follow data $y$ is found
simply from a single evaluation of the random 
(with respect to initial condition) function $\cG_\eps$.

Data is generated for the parameter set $(\sigma, r^\dagger, b^\dagger)  = (10, 28, \frac{8}{3})$, for which system \eqref{eq:lorenz63} exhibits chaotic behavior. Matrix $\Delta_{obs}$ is set to zero. 
We estimate $\Delta(\theta)$, with a matrix $\Delta_{model}$ computing a single long
trajectory of \eqref{eq:lorenz63} at $\theta^\dagger = (r^\dagger,b^\dagger)$, over $360$ time units. This is split into windows of size $10$ (neglecting the first $30$ units) and we set $\Delta_{model}$ to be the empirical covariance of $\mathcal{G}_\eps(\theta^\dagger)$ over the windows.  

We then solve the inverse problem using \od{the} time
discretization of the EKS SDE \eqref{eq:eksN}.
To ensure that there is minimal correlation between subsequent evaluations of the forward map, we set the initial condition of \eqref{eq:lorenz63}, at each step of the sampling
algorithm and for
each ensemble member, to be the state of the dynamical system from the previous ODE solve, for the same ensemble member, evaluated at a large random time $t \gg 10$.

\begin{figure}[ht]
    \centering
    \includegraphics[width=0.75\textwidth]{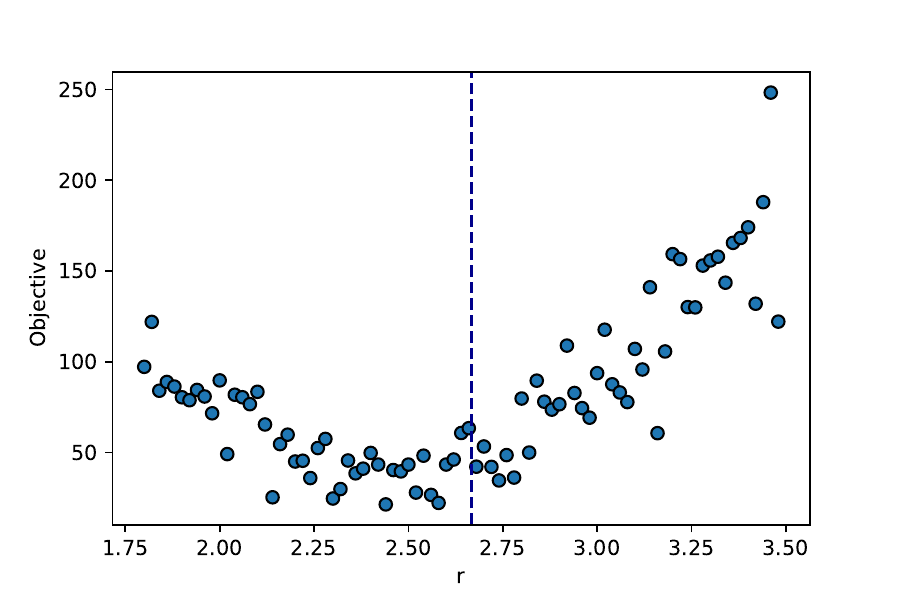}
    \caption{Profile of the noisy negative log-likelihood  over $r$ for  $b$ fixed at optimal value.  The blue dashed line denotes the 'true' value $r=8/3$. }
    \label{fig:lorenz_1d_snapshot}
\end{figure}

Given observation data $y$ and noisy forward
model $\cG_\eps$ we define the 
negative log-likelihood function  \begin{align}
    V_L(\theta) := \frac{1}{2} \langle (y-\cG_\eps(\theta)), \Delta_{model}^{-1} (y - \cG_\eps(\theta)) \rangle.
\end{align}
Figure \ref{fig:lorenz_1d_snapshot} shows the profile of $V_L$ versus $r$, for a fixed (truth) value of $b=28$. 
 We denote by $V(\theta)$ 
the negative log-posterior density found by adding
the prior quadratic form to $V_L.$
 
\begin{figure}
    \centering
    \begin{subfigure}[b]{0.45\textwidth}
    \includegraphics[width=\textwidth]{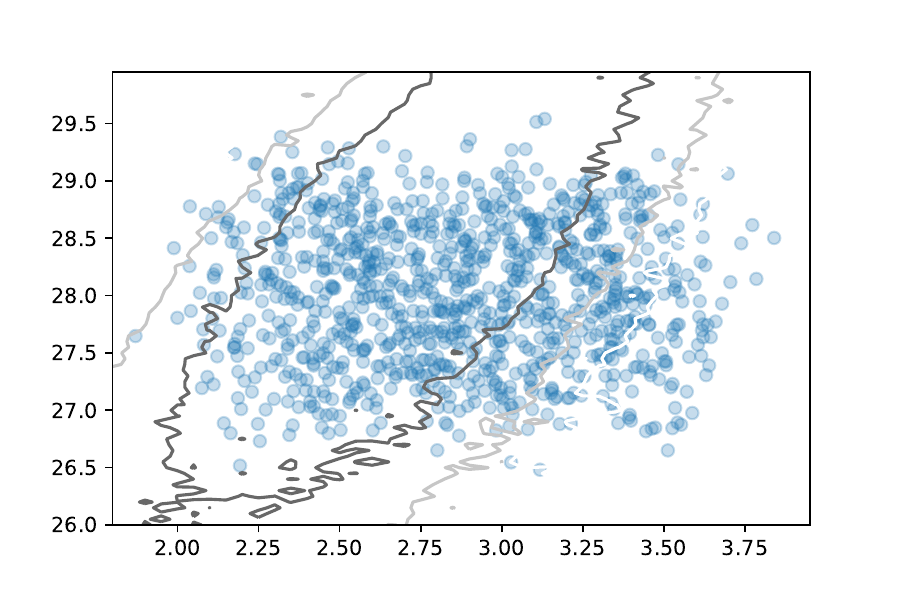}
    \caption{Ensemble Langevin Sampler}
    \end{subfigure}
    \hfill
    \begin{subfigure}[b]{0.45\textwidth}
    \includegraphics[width=\textwidth]{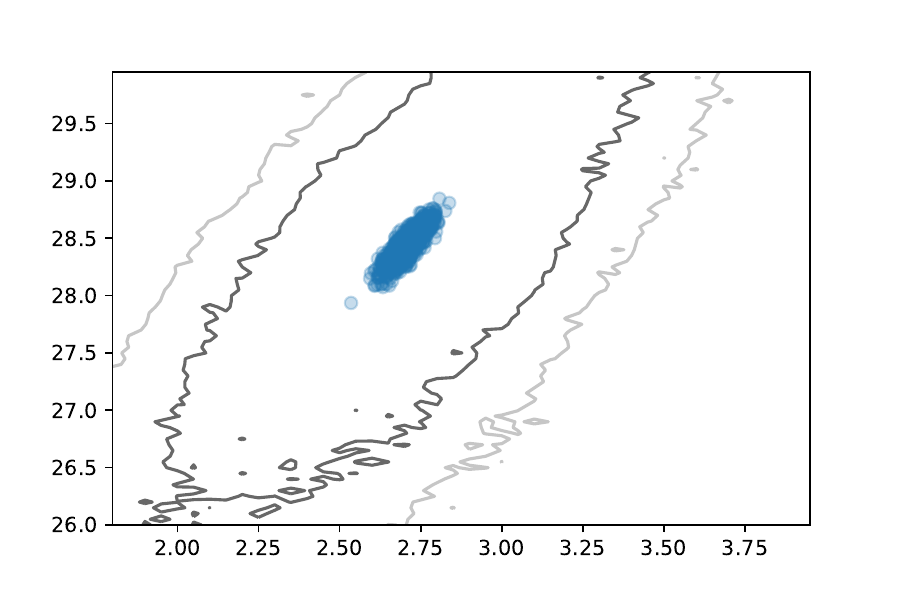}
    \caption{Ensemble Kalman Sampler}
    \end{subfigure}
    
    \begin{subfigure}[b]{0.45\textwidth}
    \includegraphics[width=\textwidth]{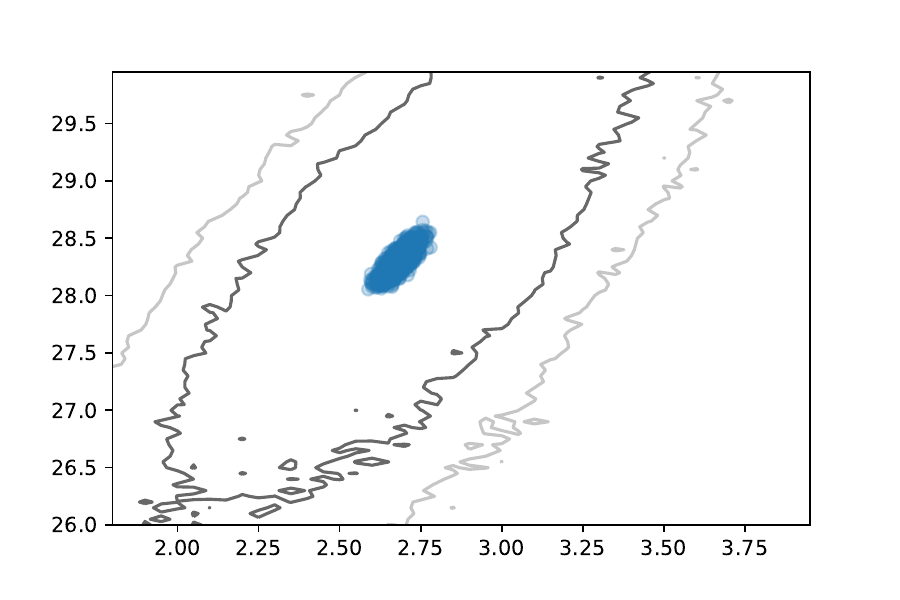}
    \caption{Ensemble GP Sampler}
    \end{subfigure}
    \hfill
    \begin{subfigure}[b]{0.45\textwidth}
    \includegraphics[width=\textwidth]{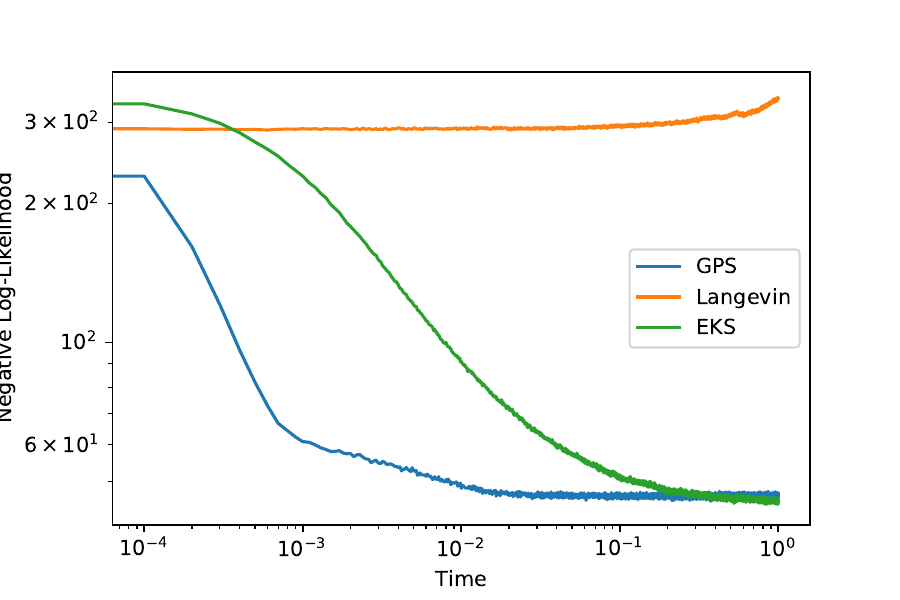}
    \caption{Negative Log-Likelihood}
    \end{subfigure}
    
    \caption{The true parameter value is $\theta=(r,b)=(28,8/3).$ Comparison of ELS, EKS
    and EGPS after simulating the ensemble for $1$ unit of time. The contour plot indicates the posterior distribution $V(\theta)$ while the dots denote the ensemble at the final time.  The bottom \as{right} plot shows the evolution of negative log likelihood as a function of time.}
    \label{fig:lorenz_hist}
\end{figure}

The EKS, ELS and EGPS processes were all simulated for one algorithmic time unit, with the step-size adjusted to ensure process stability. Each process is simulated with $N=10^3$ particles, with initial condition distributed as $U([27,29]\times[2.25,3.5])$.  In Figure \ref{fig:lorenz_hist}a)-c) we plot the particle ensemble at the final time for each method, overlaid with a contour plot of the negative log-posterior density $V(\theta)$.   In Figure \ref{fig:lorenz_hist}d) we plot the negative log-likelihood function averaged over the finite-time ensemble for each process, as the algorithms progress.   It is clear from the plots in Figures \ref{fig:lorenz_hist} that the EKS and EGPS have concentrated around the true value and are distributed according to a smoothed version of the posterior.  On the other hand, the particles undergoing ELS dynamics remain trapped around the local minima of the multiscale posterior distribution, preventing the particles from concentrating in a \as{similar fashion;
indeed the ELS is visibly close to the initial
(uniform on a rectangle) distribution of the
ensemble.

In summary the results of this subsection, where
the forward model is random, closely
mirror those from the previous subsection, where the
forward model is periodic. This substantiates our
claim that the analysis of the periodic case, contained in Sections \ref{sec:EKS} and \ref{sec:LAN}
is informative beyond the confines of the theory.}

\subsection{Multi-Modal Posteriors}
\label{ssec:multimodal}
\begin{figure}
    \centering
    \begin{subfigure}[b]{0.3\textwidth}
\includegraphics[width=\textwidth]{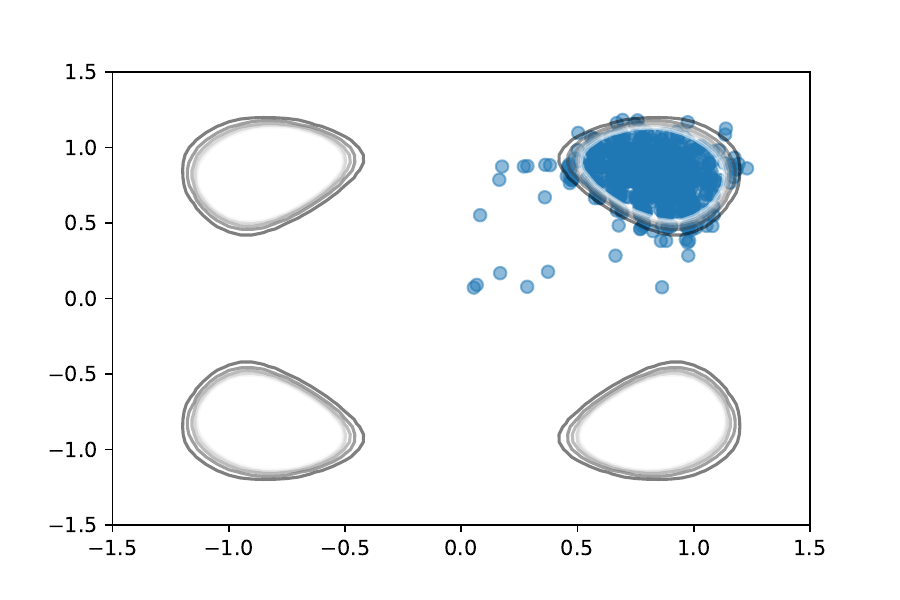}
    \caption{ELS}
    \end{subfigure}
    \hfill
    \begin{subfigure}[b]{0.3\textwidth}
    \includegraphics[width=\textwidth]{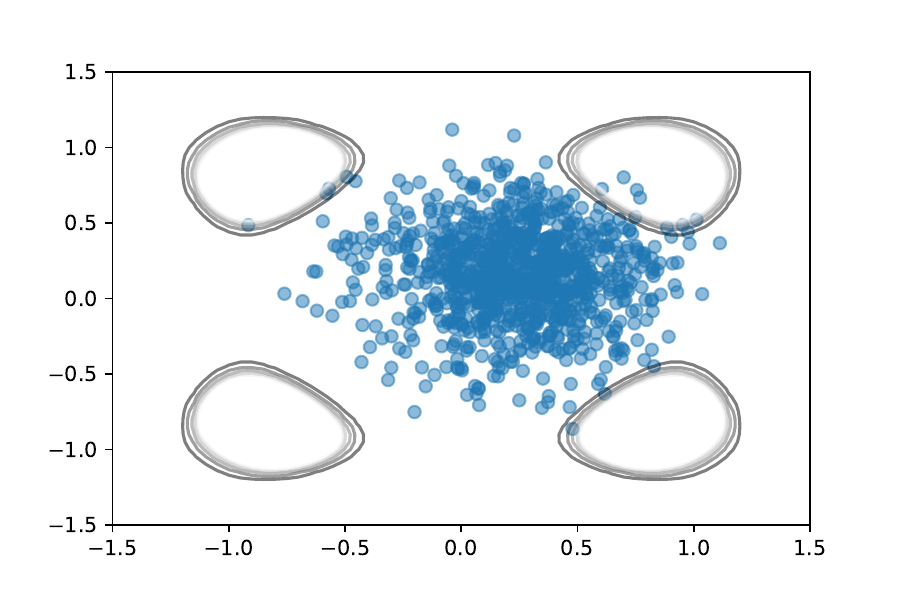}
    \caption{EKS}
    \end{subfigure}
    \hfill
    \begin{subfigure}[b]{0.3\textwidth}
    \includegraphics[width=\textwidth]{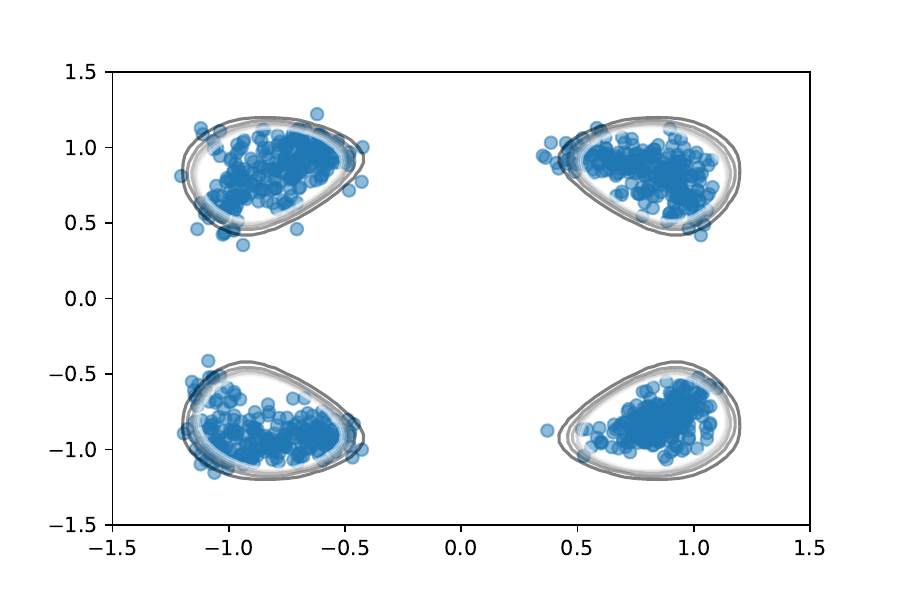}
    \caption{EGPS}
    \end{subfigure}

    \caption{Comparison of the three approaches after simulating the ensemble for $10$ units of time for the Bayesian inverse problem highlighted in Section \ref{ssec:multimodal}.  The contour plot indicates the posterior distribution associated with the slowly-varying forward map $G_0$,  while the points denote the ensemble at the final time. }
    \label{fig:multimodal}
\end{figure}

\as{It is well-known that multi-modal posteriors pose a significant challenge for ensemble Kalman-based 
approaches, since 
such approaches are constructed using a
Gaussian ansatz. Recent work on ensemble-based
interacting particle systems in
\cite{reich2019fokker} has shown the potential
for designing new interacting particle systems
which address this multi-model challenge; these methods are based on approximating nonlinear
Fokker-Planck equations arising from mean-field
dynamics, by means of particle based RKHS methods.
In the following we illustrate that, unsurprisingly,
the EGPS can achieve similar success, since it is
follows similar principles to those underlying the
work in \cite{reich2019fokker}.

We consider} the  inverse problem for the form \eqref{eq:ob0} for the unknown parameter $x\in\mathbb{R}^2$ given a multiscale forward map of the \as{form $G_\epsilon$ defined, for $x=(x_1,x_2)$, by
\begin{subequations}
\begin{align}
G_\epsilon(x) &= G_0(x) + G_1(x/\epsilon),\\ 
    G_0(x) &= (x_1^2 - 1)^2  + (x_2^2 -1)^2, \quad 
    G_1(x)= \nu (\sin(2\pi x_1)+\sin(2\pi x_2)),
\end{align}
\end{subequations}
and where $\Gamma = \gamma^2 I$.}
To demonstrate the three proposed methods we generate observation data $y \in\mathbb{R}$ for the truth $x^\dagger = (+1, -1)$, where $\nu = \frac{1}{10}$ and $\gamma^2 = 0.05$. We impose a Gaussian  $\mathcal{N}(0,\sigma^2 I)$ prior on the unknown parameter $x$, where $\sigma^2 = 0.1$.  As the slowly-varying component of the forward map is the non-injective function $G_0(x)= (x_1^2-1)^2 + (x_2^2-1)^2$, the associated posterior density exhibits $4$ global modes.  The ELS, EKS and EGPS were each simulated for an ensemble of $N=1000$ particles for $10$ time-units starting from a $U([-2,2]\times[-2,2])$ distribution.  Note that a significantly smaller step-size was selected for the Langevin sampler to ensure stability of the process.  We plot the final ensemble in Figure \ref{fig:multimodal}.   \as{As in the two previous
subsections}, we observe that the ELS struggles to explore the large-scale features of the posterior, in this case remaining concentrated on a single mode.  The effect of the multiscale perturbations can be clearly seen in the  final-time ensemble as the particle distribution appears 'corrugated' due to the influence of the sinusoidal component of the forward model. The EKS appears to be unaffected by the fine-scale structure in the forward model, but concentrates in a region at the centre of the posterior,  reflecting the fact that the EKS 
\as{is based on a Gaussian ansatz, tending to promote
unimodal distributions. Note, however, that with different initializations the EKS may concentrate on any} one of the four modes of the posterior, rather than a compromise between all of them. Finally, we observe that the EGPS sampler manages to effectively explore the large-scale structure of the posterior,  sampling from all four modes of the distribution.

\section{Conclusions}
\label{sec:CONC}
In this paper we discussed and analysed different ensemble methods for solving inverse problems with noisy and expensive likelihoods. \od{Such likelihoods commonly appear in practice, for example when using time-averaged statistics as data from a chaotic dynamical system.} A formal multiscale analysis approach was employed to characterise the influence of \as{rapid fluctuations on sampling when the objective is to explore the large-scale smoothly varying structure} of the posterior distribution. Within this framework we contrasted the long-term behaviour between sampling schemes which use gradient information and those which are gradient free, using the Ensemble Langevin Sampler (ELS) and Ensemble Kalman Sampler (EKS) as specific examples.

Both the formal analysis and computational experiments, \od{(which include both small-scale toy problems for comparison of methods, and a large-scale practical problem from climate science)} illustrate the robustness of EKS to \od{noisy and periodic} perturbations of the forward model and demonstrate its ability to efficiently characterise the underlying large scale structures of the resulting noisy posterior. This is not the case for Langevin methods, whose long time behavior is significantly impacted by the \as{rapid fluctations: these methods do not identify the 
correct smoohtly varying large-scale structure in
statistical equilibrium}, and are also slowed down by the presence of small-scale structure, \as{tending to get stuck near to the intialization of the ensemble.}. Motivated by the success of the EKS in this setting, we propose a new class of ensemble based methods, the so called ensemble Gaussian process samplers (EGPS) which are also robust to noisy perturbations of the forward model, but still employ gradient information to effectively explore the posterior distribution, and without making any assumptions on the distribution of the posterior.

While computational experiments have demonstrated the strong performance of the EGPS, it is evident that this method requires careful tuning of hyper-parameters, which is currently achieved using a preliminary tuning stage. Gaining an understanding of how to select these  parameters based on the multiscale structure of the forward map will be important for further algorithmic development.
\as{Furthermore, issues of efficiency, relating to
the frequency, in algorithmic time, with which the
Gaussian process is updated, needed to be fully
explored.} \od{Another potential `best of both worlds' solution, is to directly emulate $\mathcal{G}_\eps$ with Gaussian process and apply EKS/ELS (a philosophy taken in \cite{cleary2020calibrate,dunbar2020calibration}); direct emulation can achieve greater emulator accuracy but at an increased computational cost when $\mathcal{G}_\eps$ have a high-dimensional output space. Such algorithmic trade-offs should be investigated in different practical problems.}

On the theoretical front, it would be of interest to make the presented formal multi-scale arguments rigorous.  This might prove challenging as it would require bounds on the solution of the {\it cell problem}, a Poisson PDE which characterises the large-scale influence of small-scale perturbations.  Any such analysis would require tight lower bounds on the eigenvalues of the empirical covariance process uniformly over time. While this has been shown to be positive definite in \cite{garbuno2020affine}, obtaining quantitative lower bounds on the eigenvalues remains an open problem for future study.  Another interesting problem is to characterise the long-time behaviour of the ensemble Gaussian Process Sampler.  In particular, identifying conditions for stability and ergodicity along with quantifying the asymptotic bias are questions which we leave for future study.


\bibliographystyle{siamplain}

\bibliography{refs}

\appendix

\section{Details of the General Circulation Model}
\label{A:details_gcm}
In this section we provide explicit details of the model formulated in Section \ref{ssec:gcm}.  It is physically natural that we choose the relative humidity $RH \in [0,1]$ and relaxation time $\tau \in [0,\infty).$ In order to
accommodate Gaussian priors we introduce the transformation
\[
\theta=\mathcal{T}((\mathrm{RH}, \tau)) = \left(\mathrm{logit}(\mathrm{RH}),~\ln\left(\frac{\tau}{1~\mathrm{s}}\right)\right),
\]
which maps $[0,1] \times [0,\infty)$ into $\R^2$. 
The GCM is a single-valued function of $(RH, \tau)$, and hence of $\theta$, since $\mathcal{T}$ is invertible. We impose
Gaussian priors on $\theta\sim\mathcal{N}([0,10.17]^T,I)$, resulting in independent priors on the physical parameters $\mathcal{T}^{-1}(\theta)$ which are the logit-normal and lognormal distributions, ${\mathrm{RH}}\sim \mathrm{Logitnormal}(0, 1)$ and ${\tau}\sim \mathrm{Lognormal}(12~\mathrm{h},  (12~\mathrm{h})^2)$, respectively.

Given this, we now detail how the specific instance of data, $y$, and the
covariances $\Delta_{obs}$ and $\Delta_{model}$, are constructed.   The matrix
$\Delta_{model}$ is constructed as follows, noting that
for the atmosphere it is known that $T>15$ days is sufficient
to obtain statistical equilibrium \cite{BuiCheEmaMagSunZha19}. The data is generated from a control simulation, where the parameters are fixed at reference values, collected in the vector $\theta^\dagger$ such that $\mathcal{T}^{-1}(\theta^\dagger) = (\mathrm{RH}^\dagger = 0.7, \tau^\dagger = 2~\mathrm{h})$. 
Following \cite[Section 5]{cleary2020calibrate} we estimate the covariance $\Delta_{model}$ using long-time series data. 
To average in time, the control simulation outputs data in $1/4$-day time-steps that are then averaged over $30$-day windows to form each data sample; we generate $650$ of these samples, discarding the first $50$ to remove out-of-equilibrium initial-condition bias. We construct  $\Delta_{model}$, the variance of $\xi_{model}$,
empirically from the resulting $600$ samples.
We choose the variance $\Delta_{obs}$ of $\xi_{obs}$ as detailed in
\cite{dunbar2020calibration}; it is designed to be no more than $10\%$ of the variance arising from finite time-averaging,
and also to ensure physically reasonable data $y$ (e.g. precipitation data $\geq 0$), with high probability.
The data $y$ is then constructed
by drawing at random one of the $600$ $30$-day samples, representing a draw from $\mathcal{G}_0(\theta^\dagger)+\xi_{model}$, and adding to it a draw $\xi_{obs}\sim \mathcal{N}(0,\Delta_{obs})$, representing 
observation error. 
In particular we have unified $\epsilon^{-1}=30$ days.

\section{Multiscale Analysis For EKS}
\label{A:1}

In this section we derive Formal Perturbation Result \ref{fpr:1} concerning
averaging for the mean field limit of the EKS. To carry out the
analysis we extend the spatial domain of the mean field system from $\Rd$ to $\Rd \times \Td$ as is standard in the perturbation
approach described in \cite{bensoussan2011asymptotic,pavliotis2008multiscale}.
The analysis will be stream-lined
by making the following definitions. For economy of notation
we reuse the notation
$\rho(\cdot,t), \CF(\cdot,\rho)$, now to denote functions with input domain extended
from $\Rd$ to $\Rd \times \Td$; specifically, this naturally generalizes the definitions of $\rho(\cdot,t), \CF(\cdot,\rho)$
in Section \ref{sec:EKS}. 

In the following $\pi: \Rd \times \Td  \to \R^+$ denotes a probability density on $\Rd \times \Td$ and $\pi_0: \Rd  \to \R^+$ denotes a probability density on $\Rd$. Using
this notation we define
\begin{align*}
\bx(\pi)&=\intrt x\pi(x,z)dx dz,\\
\bG_0(\pi)&=\intrt G_0(x) \pi(x,z)dx dz,\quad
\bG_1(\pi)=\intrt G_1(z) \pi(x,z)dx dz,\\
\CC(\pi)&=\intrt \bigl(x-\bx(\pi)\bigr)\otimes\bigl(x-\bx(\pi)\bigr) \pi(x,z) dx dz,\\
\CF(x,z,\pi)&=\intrt \langle G_0(x')+G_1(z')-\bG_0(\pi)-\bG_1(\pi), G_0(x)+G_1(z)-y \rangle_{\Gamma} x'\pi(x',z')dx' dz',\\
\CC_0(\pi_0)&=\intr \bigl(x-\bx(\pi_0)\bigr)\otimes\bigl(x-\bx(\pi_0)\bigr) \pi_0(x) dx,\\
\CF_0(x,\pi_0)&=\intr \langle G_0(x')-\bG_0(\pi_0), G_0(x)-y \rangle_{\Gamma} x'\pi_0(x')dx'.
\end{align*}
Note that in employing this notation $\CC$, viewed as a matrix-valued functional on densities, is extended
from its definition in Section \ref{sec:EKS} to now act on densities on $\Rd \times \Td$. However if density $\rho$
is constant in $z$ then we recover the definition of $\CC(\rho)$
from Section \ref{sec:EKS}; in this case $\CC(\rho)=\CC_0(\rho).$
We also define the following differential operators.
\begin{align*}
    \Bz(\rho)\bullet&=\nabla_z \cdot\bigl(\nabla_z \cdot \bigl(\CC(\rho)\bullet\bigr)\bigr),\\
    \Bo(\rho)\bullet&=2\nabla_z \cdot\bigl(\nabla_x \cdot \bigl(\CC(\rho)\bullet\bigr)\bigr)+\nabla_z \cdot \bigl(\CF(x,z,\rho)\bullet\bigr),\\
    \Bt(\rho)\bullet&=\nabla_x \cdot\bigl(\nabla_x \cdot \bigl(\CC(\rho)\bullet\bigr)\bigr)+\nabla_x \cdot \bigl(\CF(x,z,\rho)\bullet\bigr),
\end{align*}
Under the Assumption \ref{a:1} on $G=G_{\eps}$ the finite
particle system \eqref{eq:eksN} is
\begin{align}
\label{eq:eksNx}
dX_t^i = -\Bigl(\frac{1}{N}\sum_{n=1}^{N}\langle G_0(X_t^n)+G_1(X_t^n/\eps) - \overline{G}_{0,t}-\overline{G}_{1,t}, G_0(X_t^i)+
G_1(X_t^i/\eps)- y\rangle_{\Gamma} X_t^n \Bigr)\,dt&\nonumber\\
- C_t\Sigma^{-1}X_t^i \,dt + \frac{d+1}{N}(X_t^i - \overline{X}_t)\, dt + \sqrt{2C_t}\,dW^i_t&,
\end{align}
where
\begin{align*}
 \overline{X}_t = \frac{1}{N}\sum_{n=1}^N X_t^n, \quad
 \overline{G}_{0,t} = \frac{1}{N}\sum_{n=1}^{N} G_0(X_t^n), \quad
 \overline{G}_{1,t} = \frac{1}{N}\sum_{n=1}^{N} G_1(X_t^n/\eps),
 \end{align*}
 and
 $C_t = \frac{1}{N}\sum_{n=1}^{N}\left(X_t^n - \overline{X}_t\right)\otimes \left(X_t^n - \overline{X}_t\right).$
If we introduce $Z_t^i=X_t^i/\eps$ then we obtain
\begin{subequations}
\label{eq:eksNz}
\begin{align}
dX_t^i &= -\Bigl(\frac{1}{N}\sum_{n=1}^{N}\langle G_0(X_t^n)+G_1(Z_t^n) - \overline{G}_{0,t}-\overline{G}_{1,t}, G_0(X_t^i)+
G_1(Z_t^i)- y\rangle_{\Gamma} X_t^n \Bigr)\,dt\nonumber \\
&\quad\quad\quad\quad\quad\quad\quad\quad\quad\hspace{0.5in}- C_t\Sigma^{-1}X_t^i \,dt + \frac{d+1}{N}(X_t^i - \overline{X}_t)\, dt + \sqrt{2C_t}\,dW^i_t,\\
\eps dZ_t^i &= -\Bigl(\frac{1}{N}\sum_{n=1}^{N}\langle G_0(X_t^n)+G_1(Z_t^n) - \overline{G}_{0,t}-\overline{G}_{1,t}, G_0(X_t^i)+
G_1(Z_t^i)- y\rangle_{\Gamma} X_t^n \Bigr)\,dt\nonumber \\
&\quad\quad\quad\quad\quad\quad\quad\quad\quad\hspace{0.5in}- C_t\Sigma^{-1}X_t^i \,dt + \frac{d+1}{N}(X_t^i - \overline{X}_t)\, dt + \sqrt{2C_t}\,dW^i_t,
\end{align}
\end{subequations}
where now we may write
\begin{align*}
 \overline{G}_{1,t} = \frac{1}{N}\sum_{n=1}^{N} G(Z_t^n).
 \end{align*}

Now consider the mean field SDE defined by this system.
Similarly to the exposition in Section \ref{sec:EKS} this takes
the form
\begin{subequations}
\label{eq:FS}
\begin{align}
    dx&=-\CF(x,z,\rho)dt +\sqrt{2\CC(\rho)}dW,\\
    \eps dz&=-\CF(x,z,\rho)dt +\sqrt{2\CC(\rho)}dW,
\end{align}
\end{subequations}
where, to be self-consistent the density $\rho(x,z,t)$ must satisfy the equation
\begin{equation}
    \label{eq:MSFP}
    \partial_t \rho=\frac{1}{\eps^2}\Bz(\rho)\rho+\frac{1}{\eps}\Bo(\rho)\rho+
    \Bt(\rho)\rho.
\end{equation}

We seek a solution in the form
\begin{equation}
    \label{eq:exp}
    \rho=\rho_0+\eps \rho_1+\eps^2 \rho_2 + \cdots
\end{equation}
and assume the normalizations
\begin{subequations}
    \label{eq:norm}
    \begin{align*}
    \intrt \rho_0(x,z,t)dxdz&=1,\\
    \intrt \rho_j(x,z,t)dxdz&=0, \quad j \ge 1;
    \end{align*}
\end{subequations}
this ensures that $\rho$ integrates to $1$.
We now expand the operators $\Bz(\rho), \Bo(\rho), \Bt(\rho)$ about
$\rho_0.$ To this end we first note that
\begin{equation}
    \label{eq:Cexp}
    \CC(\rho)=\CC(\rho_0)+\eps \, \mathcal{C}_1+\eps^2 \,\mathcal{C}_2;
\end{equation}
we will not need the precise forms of $\mathcal{C}_1$ and $\mathcal{C}_2$ in what follows.
From this we deduce that
\begin{subequations}
\label{eq:bbz}
\begin{align}
    \Bz(\rho)\bullet&=\Bz(\rho_0)\bullet+\eps \nabla_z\cdot\bigl(\nabla_z \cdot( \mathcal{C}_1\bullet)\bigr)+
    \eps^2\nabla_z\cdot\bigl(\nabla_z \cdot( \CC_2\bullet)\bigr),\\
    \Bo(\rho)\bullet&=\Bo(\rho_0)\bullet+2\eps \nabla_z\cdot\bigl(\nabla_x \cdot( \CC_1\bullet)\bigr)+
    \eps\nabla_z\cdot\bigl(D_{\rho}\CF(x,z,\rho_0)\rho_1\bullet\bigr).
\end{align}
\end{subequations}
Using these expressions, and substituting \eqref{eq:exp} into 
\eqref{eq:MSFP} and equating terms of size $\mathcal{O}(\eps^{-2})$, $\mathcal{O}(\eps^{-1})$ and $\mathcal{O}(1)$ respectively, gives the following equations:
\begin{subequations}
\label{eq:FA}
\begin{align}
\Bz(\rho_0) \rho_0&=0,\\
\Bz(\rho_0) \rho_1&=-\Bo(\rho_0) \rho_0 -\nabla_z \cdot\bigl(\nabla_z \cdot(\CC_1 \rho_0)\bigr),\\
\begin{split}
\Bz(\rho_0) \rho_2&=-\Bo(\rho_0) \rho_1 -\nabla_z \cdot\bigl(\nabla_z \cdot(\CC_2 \rho_0)\bigr)  -2\nabla_z \cdot\bigl(\nabla_z \cdot(\CC_1 \rho_1)\bigr) \\
&\phantom{=-Bo(\rho_0) \rho_1 }
-\nabla_z \cdot\bigl(D_{\rho}\CF(u,v,\rho_0)\rho_1\bigr)-\Bt(\rho_0)\rho_0+\partial_t \rho_0.
\end{split}
\end{align}
\end{subequations}
Note that $\Bz(\rho_0)$ is a differential operator in $z$ only
and that its nullspace comprises
constants in $z$. We see that equation (\ref{eq:FA}a) is solved by assuming that $\rho_0(x,t)$
only, and is independent of $z$, because
$\Bz(\rho_0)$ has constants with respect to $z$ in its null-space. We now turn to equation (\ref{eq:FA}b), noting that the operator
$\Bz(\rho_0)$ is self-adjoint. Thus the Fredholm alternative requires that
the right-hand side of equation (\ref{eq:FA}b) is orthogonal to constants
on $\Td$ in $z$ for a solution $\rho_1$ to exist; this is a condition which is automatically satisfied because the right hand side is a divergence
with respect to $z$. Using this structure we find a solution $\rho_1$ which we make unique by imposing (\ref{eq:norm}b).
We again apply the Fredholm alternative, now to ensure
existence of a solution of equation (\ref{eq:FA}c). The condition that the right hand
side is orthogonal to constants
on $\Td$ in $z$ then gives, noting that divergences in $z$ again contribute nothing,
$$\partial_t \rho_0=\nabla_x \cdot \nabla_x\cdot\Bigl(\CC_0(\rho_0)\rho_0\Bigr)+\nabla_x\cdot\Bigl(\intt \CF(x,z,\rho_0)dz
\rho_0\Bigr).$$
Using the fact that $\rho_0$ is independent of $z$, and since
$G_1$ has mean zero on $\Td$, it follows that
$$\partial_t \rho_0=\nabla_x \cdot \nabla_x\cdot\Bigl(\CC_0(\rho_0)\rho_0\Bigr)+\nabla_x\cdot\Bigl( \CF_0(x,\rho_0)
\rho_0\Bigr).$$
This is the nonlinear Fokker-Planck equation \eqref{eq:NLFP0} associated with the desired
averaged mean-field limit equations, after noting that $\CC_0(\rho_0)$ is the same as $\CC(\rho_0)$, with the latter
using the notation for matrix-valued functional $\CC$ as defined
in Section \ref{sec:EKS}.

\section{Multiscale Analysis For Ensemble Langevin Dynamics}
\label{A:2}
In this section we derive Formal Perturbation Result 2. This result concerns homogenization of the mean field limit for ensembles of coupled particles undergoing overdamped Langevin dynamics defined by noisy
forward model $G_\epsilon.$ In the
mean field limit the particle is ergodic with respect to $\pi\propto e^{-V_\epsilon}$, where 
$$
  V_{\epsilon}(x) :=  \frac{1}{2}\langle (y - G_\epsilon(x)), \Gamma^{-1}(y - G_{\epsilon}(x))\rangle + \frac{1}{2}\langle x, \Sigma^{-1} x \rangle,
$$
and $G_\epsilon$ is given by \eqref{eq:forward}. The mean-field density $\rho$ satisfies the following nonlinear PDE
\begin{equation}
\label{eq:mfl}
    \partial_t \rho = \nabla_x\cdot\left(\mathcal{M}(\rho)\left(\nabla_x \rho + \nabla_x V_{\epsilon} \rho\right)\right),
\end{equation}
where $\mathcal{M}$ is a bounded linear operator on the vector-valued Hilbert space  $L^2(\mathbb{R}^d; \mathbb{R}^d)$. 
\\\\ 
We write $V_{\epsilon}(x) = V(x,x/\epsilon) = V_0(x) + V_1(x,x/\epsilon)$ where 
$$
  V_{0}(x)=\frac{1}{2}\langle (y - G_0(x)), \Gamma^{-1}(y - G_{0}(x))\rangle + \frac{1}{2}\langle x, \Sigma^{-1} x \rangle,
$$
and 
$$
    V_1(x,z) = \frac{1}{2}\langle G_1(x,z), \Gamma^{-1}G_{1}(x,z)\rangle + \langle (y - G_0(x)), \Gamma^{-1}G_1(x,z)\rangle.
$$
Also writing $\nabla_x \mapsto
\nabla_x+\epsilon^{-1}\nabla_z$ in
\eqref{eq:mfl}, and viewing $\rho$ as a function of $(x,z,t)$, we can rewrite the
nonlinear Fokker-Planck equation as
\begin{equation}
\label{eq:first}
    \partial_t \rho = \frac{1}{\epsilon^2}\mathcal{B}_0(\rho)\rho +  \frac{1}{\epsilon}\mathcal{B}_1(\rho)\rho +  \mathcal{B}_2(\rho)\rho,
\end{equation}
where 
\begin{align*}
\mathcal{B}_0(\rho)\bullet&= \nabla_z\cdot\left(\mathcal{M}(\rho)\left(\nabla_z \bullet+ \nabla_z V\bullet\right)\right) \\
\mathcal{B}_1(\rho)\bullet&= \nabla_x\cdot\left(\mathcal{M}(\rho)\left(\nabla_z \bullet+ \nabla_z V\bullet\right)\right) + \nabla_z\cdot\left(\mathcal{M}(\rho)\left(\nabla_x \bullet+ \nabla_x  V\bullet\right)\right)\\
\mathcal{B}_2(\rho)\bullet&=  \nabla_x\cdot\left(\mathcal{M}(\rho)\left(\nabla_x \bullet+ \nabla_x V\bullet\right)\right).
\end{align*}
As in Appendix \ref{A:1} we have extended the spatial domain of the mean field equation from $\mathbb{R}^d$ to $\mathbb{R}^d\times \mathbb{T}^d$ and $\rho(\cdot,\cdot,t)$ is a probability density 
function on $\mathbb{R}^d\times \mathbb{T}^d$, for each fixed $t$.  

Similarly to the analysis in Appendix \ref{A:1}, 
$\mathcal{B}_0(\rho)$ is a differential operator in $z$ only, but now the null-space has non-trivial variation in $z$: it
comprises functions of the form
$\exp\bigl(-V_1(x,z)\bigr).$ In this homogenization setting we should not expect the leading order term of the solution, $\rho_0$ to be independent of the fast scale fluctuations, nor should we expect pointwise convergence of $\rho$ to $\rho_0$.  We thus introduce the following rescaling of the standard perturbation expansion to account for the fast-scale fluctuations in $\rho$, as in
\cite{givon2004extracting}[Section 6.2]:
\begin{subequations}
    \label{eq:second}
    \begin{align}
\rho &= \rho_0+  \epsilon  \rho_1  +  \epsilon^2   \rho_2 +  \ldots\ \\
 &=  e^{-V}\bigl( \chi_0+  \epsilon  \chi_1  +  \epsilon^2   \chi_2 +  \ldots\bigr)
 \end{align}
\end{subequations}
where $\chi_i=\chi_i(x,z,t)$ and $V(x,z)=V_0(x)+V_1(x,z).$ We
impose the conditions 
\begin{align*}
    \int_{\mathbb{R}^d}\int_{\mathbb{T}^d} \chi_0(t,x,z)e^{-V(x,z)}\,dx \,dz&= 1,  \\
    \int_{\mathbb{R}^d}\int_{\mathbb{T}^d} \chi_j(t,x,z)e^{-V(x,z)}\,dx &= 0, \quad j \geq 1.
\end{align*}
We have $\rho_0(x,z,t)=e^{-V(x,z)}\chi_0(x,z,t).$ Similarly to the derivation in Section \ref{A:1},
we assume that $\mathcal{M}$ admits the following regular expansion:
\begin{equation} 
    \mathcal{M}(\rho)= \mathcal{M}(\rho_0)+\eps \mathcal{M}_1 + \eps^2 \mathcal{M}_2 + \ldots
\end{equation}
where $\mathcal{M}_1$ and $\mathcal{M}_2$ are independent of $\epsilon$.   In particular, both the possible choices of $\mathcal{M}$ identified in Section \ref{sec:LAN} admit such an expansion. From this we observe that we can express $\mathcal{B}_0(\rho)\bullet$ and $\mathcal{B}_1(\rho)\bullet$ in terms of $\mathcal{B}_0(\rho_0)\bullet$ and $\mathcal{B}_1(\rho_0)\bullet$ respectively, as follows:
\begin{align*}
    \mathcal{B}_0(\rho) \bullet&= \mathcal{B}_0(\rho_0)\bullet+ \epsilon \mathcal{B}_0^{(1)}\bullet+ \epsilon^2 \mathcal{B}_0^{(2)} \bullet+\ldots\\
    \mathcal{B}_1(\rho)\bullet&= \mathcal{B}_1(\rho_0)\bullet+  \epsilon \mathcal{B}_1^{(1)} \bullet+ \dots\\
    \mathcal{B}_2(\rho)\bullet&= \mathcal{B}_2(\rho_0)\bullet+ \dots
\end{align*}
where the linear operators $\{\mathcal{B}_i^{(j)}\}$ acting on the space of probability density functions
are defined by
\begin{align*}
    \mathcal{B}_0^{(1)}\bullet&= 
     \nabla_z\cdot\left(\mathcal{M}_1(\nabla_z \bullet + \nabla_z V \bullet))\right) \\   
     \mathcal{B}_0^{(2)}\bullet&=  \nabla_z\cdot\left(\mathcal{M}_2(\nabla_z \bullet + \nabla_z V \bullet))\right) \\
    \mathcal{B}_1^{(1)}\bullet&=  \nabla_x\cdot\left(\mathcal{M}_1\left(\nabla_z \bullet + \nabla_z V\bullet\right)\right) + \nabla_z\cdot\left(\mathcal{M}_1\left(\nabla_x \bullet + \nabla_x  V\bullet\right)\right)\\
  \mathcal{B}_1^{(2)}\bullet&=  \nabla_x\cdot\left(\mathcal{M}_2\left(\nabla_z \bullet + \nabla_z V\bullet\right)\right) + \nabla_z\cdot\left(\mathcal{M}_2\left(\nabla_x \bullet + \nabla_x  V\bullet\right)\right).\textbf{}
\end{align*}
Using these expressions in \eqref{eq:first}, substituting the expansion \eqref{eq:second} and equating terms of size $\mathcal{O}(\eps^{-2}), \mathcal{O}(\eps^{-1})$ and $\mathcal{O}(1)$ respectively, gives the following equations:
\begin{align}
    \label{eq:fa1}\mathcal{B}_0(\rho_0)\rho_0 &= 0\\
     \label{eq:fa2}\mathcal{B}_0(\rho_0)\rho_1 &= -\mathcal{B}_1(\rho_0)\rho_0 - \mathcal{B}_0^{(1)}\rho_0 \\
     \label{eq:fa3}\mathcal{B}_0(\rho_0)\rho_2 &=  \partial_t \rho_0  
     -\mathcal{B}_1(\rho_0)\rho_1 - \mathcal{B}_0^{(1)}\rho_1
    - \mathcal{B}_1^{(1)}\rho_0 - \mathcal{B}_0^{(2)}\rho_0 
     - \mathcal{B}_2(\rho_0)\rho_0.
\end{align}
Noting that $\nabla_z V=\nabla_z V_1$ it follows that the $\mathcal{O}(\epsilon^{-2})$ equation \eqref{eq:fa1} can be expressed as 
$$
\nabla_z\cdot\left(\mathcal{M}(\rho_0)e^{-V_1}\nabla_z \chi_0 \right) = 0.
$$
This equation may be solved by noting that $\chi_0$ must be a constant with respect to $z$
since the operator acting on $\chi_0$ has only constants in its nullspace; thus $\chi_0(x,z,t) = \chi_0(x,t)$. The second equation \eqref{eq:fa2} for the $\mathcal{O}(\epsilon^{-1})$ terms gives
$$
\nabla_z\cdot\left(\mathcal{M}(\rho_0)e^{-V_1}\nabla_z \chi_1 \right) = -\nabla_z\cdot\left(\mathcal{M}(\rho_0)e^{-V_1}\nabla_x \chi_0\right).
$$
The operator acting on $\chi_1$ is self-adjoint with only constants in its null-space. Thus, by the Fredholm alternative the
equation has a solution since the right hand side is divergence free. We can write this solution in the form $\chi_1 = \chi\cdot \nabla_x\chi_0$ where $\chi$ satisfies the following PDE:
$$
\nabla_z\cdot\left(\mathcal{M}(\rho_0)e^{-V_1}(\nabla_z \chi + I) \right) = 0. 
$$

Multiplying this identity by $\chi$ and integrating by parts implies the following identity
which we will use at the end of this section to study properties of the homogenized limit:
\begin{equation}
\label{eq:ident}
\intT \nabla_z \chi^{\top}  \mathcal{M}(\rho_0)\nabla_z \chi e^{-V_1}\,dz=-\intT   \mathcal{M}(\rho_0)\nabla_z \chi e^{-V_1}\,dz.
\end{equation}

We now consider the $\mathcal{O}(1)$ terms, and equation \eqref{eq:fa3}. Again invoking the Fredholm 
alternative requires that the integral of the RHS integrates to zero with respect to $z$. We note that every term
appearing in the expression
$$\mathcal{B}_0^{(1)}\rho_1
    + \mathcal{B}_1^{(1)}\rho_0 + \mathcal{B}_0^{(2)}\rho_0$$
is a divergence with respect to $z$ with the exception of one divergence with respect to $x$ which is
identically zero. It follows that
$$\int_{\mathbb{T}^d} \Bigl(\partial_t \rho_0  
     -\mathcal{B}_1(\rho_0)\rho_1
     - \mathcal{B}_2(\rho_0)\rho_0\Bigr)dz=0.$$ 
Evaluating this integral, 
we obtain
\begin{align*}
    \intT \partial_t \rho_0\,dz &= \int_{\mathbb{T}^d} \nabla_x\cdot\bigl(\mathcal{M}(\rho_0)(\nabla_z \rho_1 + \nabla_z V\rho_1)\bigr)\,dz  \\
    & + \int_{\mathbb{T}^d} \nabla_z\cdot\left(\mathcal{M}(\rho_0)\left(\nabla_x \rho_1  +  \nabla_x  V\rho_1\right)\right)\,dz \\
    & +  \int_{\mathbb{T}^d} \nabla_x\cdot\bigl(\mathcal{M}(\rho_0)(\nabla_x \rho_0 + \nabla_x V\rho_0)\bigr)\,dz.
\end{align*}
The second term on the RHS drops out by the divergence theorem.  Noting that $\mathcal{M}(\rho_0)$ does not depend on the fast variable $z$ we then obtain
\begin{align}
\label{eq:o2_1}
\begin{split}
    \int_{\mathbb{T}^d} \partial_t \rho_0\,dz &=  \nabla_x\cdot\left( \mathcal{M}(\rho_0)\int_{\mathbb{T}^d}\left(\nabla_z \rho_1  + \nabla_z V\rho_1\right)dz\right)\\
    &+   \nabla_x\cdot\left(\mathcal{M}(\rho_0)\int_{\mathbb{T}^d}\left(\nabla_x \rho_0 + \nabla_x V\rho_0\right)dz.\right).
    \end{split}
\end{align}
Substituting  $\rho_1 =\chi_1 e^{-V} = \chi\cdot\nabla_x \chi_0 e^{-V}$ we obtain
\begin{align*}
  \nabla_z \rho_1 + \nabla_z V \rho_1  = e^{-V}\nabla_z\chi_1 = e^{-V}\nabla_z \chi \nabla_x \chi_0.
\end{align*}
Similarly substituting $\rho_0 = \chi_0 e^{-V}$ yields
$$
\nabla_x\rho_0 + \nabla_x V\rho_0 = e^{-V}\nabla_x(\rho_0 e^V) = e^{-V}\nabla_x\chi_0. 
$$
Thus equation \eqref{eq:o2_1} becomes
\begin{align*}
    \partial_t \chi_0e^{-V_0} \int_{\mathbb{T}^d} e^{-V_1}\,dz &= \nabla_x\cdot\left(e^{-V_0}\left[\int_{\mathbb{T}^d} \mathcal{M}(\rho_0)e^{-V_1}\nabla_z \chi\,dz\right]\nabla_x \chi_0\right) \\
    &\quad +\nabla_x\cdot\left(e^{-V_0}\left[\int_{\mathbb{T}^d} \mathcal{M}(\rho_0)e^{-V_1}\,dz\right]\nabla_x\chi_0\right) 
\end{align*}
We now define
$$
Z(x) = \int_{\mathbb{T}^d} e^{-V_1(x,z)}\,dz,
$$
and 
$$
\CD(\rho_0,x) = \frac{1}{Z(x)} \int_{\mathbb{T}^d} \mathcal{M}(\rho_0)e^{-V_1}(I + \nabla_{z}\chi)\,dz;
$$
we write $Z$ and $\CD(\rho_0)$, suppressing explicit dependence on $x$ in some of what follows. The effective dynamics becomes
\begin{align*}
    \partial_t \chi_0 e^{-V_0}Z &= \nabla_x\cdot\bigl(e^{-V_0}Z\CD(\rho_0)\nabla_x \chi_0\bigr).
\end{align*}
To conclude the limit argument, let $\phi \in C^{\infty}_C$, then we have formally that
\begin{align*}
\int_{\mathbb{T}^d}	\int_{\mathbb{R}^d} \phi(x,t) \rho(x,t)\,dx\,dz & \xrightarrow{\epsilon\rightarrow 0} \int_{\mathbb{T}^d}\int_{\mathbb{R}^d} \phi(x)\rho_0(x,z, t)dx\,dz\\ & =  \int_{\mathbb{R}^d}\int_{\mathbb{T}^d} \phi(x)\chi_0(x,t)e^{-V(x,z)}\,dz\,dx\\ &=  \int_{\mathbb{R}^d} \phi(x)\chi_0(x,t)\int_{\mathbb{T}^d}e^{-V(x,z)}\,dz\,dx\\
&=  \int_{\mathbb{T}^d} \int_{\mathbb{R}^d} \phi(x)\chi_0(x,t)e^{-V_0(x)}Z(x) \,dx\,dz.
\end{align*}
Thus we identify the leading order term in the expansion for the density as  $\rho_0(x,t)  = \chi_0(t,x) e^{-V_0(x)}Z(x)$.  This is the average over the torus of $\rho_0(x,z,t)$ and
we overload notation for $\rho_0$ deliberately to avoid proliferation of symbols.
The equation for $\rho_0=\rho_0(x,t)$ is
$$
\partial_t  \rho_0 = \nabla_x\cdot\Bigl(\CD(\rho_0)e^{-V_0}Z\nabla_x\bigl(\rho_0/(e^{-V_0}Z)\bigr)\Bigr) =  \nabla_x\cdot\bigl(\CD(\rho_0)(\nabla_x \rho_0  + \nabla_x \overline{V}\rho_0\bigr),
$$
where $\overline{V}(x)  = V_0(x)-\log Z(x)$.  This  is the mean-field limit for a system of overdamped Langevin particles evolving in a potential  $\overline{V}$ with density dependent  diffusion tensor $\CD(\rho_0)$.
\\\\
The equilibrium solution of this equation is given by  $\overline{\pi}(dx)  \propto  e^{-V_0(x)}Z(x)$.  For general forward problems, this will  be  different from the posterior distribution $\pi_0(dx)\propto e^{-V_0(x)}\,dx$ associated with the unperturbed forward model $G_0$.  
\\\\
Furthermore we  also  note the introduction  of  a  slow-down in  the  evolution of $\rho_0(\cdot,t)$  to equilibrium  as $t\rightarrow  \infty$, in comparison with
the original ensemble Langevin dynamics in the smooth potential $V_0.$ This may be seen by comparing the linear operator $\CD(\cdot, x)$, arising in the homogenized equation for $\rho_0$ with $\mathcal{M}(\cdot)$.  Indeed, using
\eqref{eq:ident}, we can rewrite this {\it effective diffusion} operator
$\CD(\rho_0, x)$ as
\begin{align*}
    \CD(\rho_0, x) &= \frac{1}{Z(x)}\int_{\mathbb{T}^d} \mathcal{M}(\rho_0)e^{-V_1}(I + \nabla_z \chi)\,dz\\ & = \frac{1}{Z(x)}\int_{\mathbb{T}^d} \mathcal{M}(\rho_0)e^{-V_1(x,z)}\,dz - \int_{\mathbb{T}^d}\nabla_z \chi^{\top} \mathcal{M}(\rho_0)\nabla_z \chi e^{-V_1(x,z)}\,dz \\
    &= \mathcal{M}(\rho_0) - \int_{\mathbb{T}^d}\nabla_z \chi^{\top} \mathcal{M}(\rho_0)\nabla_z \chi e^{-V_1(x,z)}\,dz.
\end{align*}
Thus, for arbitrary $\zeta \in L^2(\mathbb{R}^d; \mathbb{R}^d)$, \eqref{eq:oqf} holds.
This demonstrates that the effective diffusion is always smaller than  or equal to that
in the potential defined by $G_0$, in the sense of spectrum. For a single particle
in a multiscale potential, this slowing down phenomenon is analyzed in \cite{Olla94homogenizationof}.

\end{document}